\documentclass[12pt]{amsart}
\usepackage{graphicx}
\usepackage{amsmath}
\usepackage{amsfonts}
\usepackage{amssymb}
\usepackage{mathtools}
\usepackage{amscd}
\usepackage[all,cmtip]{xy}
\usepackage{amsthm}
\DeclareFontFamily{U}{mathb}{\hyphenchar\font45}
\DeclareFontShape{U}{mathb}{m}{n}{
  <-6> mathb5 <6-7> mathb6 <7-8> mathb7
  <8-9> mathb8 <9-10> mathb9
  <10-12> mathb10 <12-> mathb12
}{}
\DeclareSymbolFont{mathb}{U}{mathb}{m}{n}
\DeclareMathSymbol{\Asterisk}{\mathbin}{mathb}{"06}
\DeclareMathSymbol{\righttoleftarrow}{\mathrel}{mathb}{"FD}
\usepackage{geometry}
\usepackage[colorlinks=true,linkcolor=blue,citecolor=blue,urlcolor=blue]{hyperref}
\usepackage[capitalise]{cleveref}
\usepackage{xcolor}
\usepackage{comment}

\newtheorem{theorem}{Theorem}[section]
\newtheorem{proposition}[theorem]{Proposition}
\newtheorem{corollary}[theorem]{Corollary}
\newtheorem{lemma}[theorem]{Lemma}

  \theoremstyle{definition}
\newtheorem{definition}[theorem]{Definition}
\newtheorem{convention}[theorem]{Convention}
\newtheorem{example}[theorem]{Example}

\newtheorem{observation}[theorem]{Observation}

\crefname{theorem}{Theorem}{Theorems}
\crefname{proposition}{Proposition}{Propositions}
\crefname{corollary}{Corollary}{Corollaries}
\crefname{lemma}{Lemma}{Lemmas}
\crefname{claim}{Claim}{Claims}
\crefname{fact}{Fact}{Facts}
\crefname{definition}{Definition}{Definitions}
\crefname{convention}{Convention}{Conventions}
\crefname{example}{Example}{Examples}
\crefname{remark}{Remark}{Remarks}
\crefname{question}{Question}{Questions}
\crefname{notation}{Notation}{Notations}
\crefname{observation}{Observation}{Observations}

\newcommand{\dbN}{\mathbb{N}}

\newcommand{\dbR}{\mathbb{R}}

\newcommand{\dbZ}{\mathbb{Z}}

\newcommand{\calC}{{\mathcal C}}

\newcommand{\calS}{{\mathcal S}}
\newcommand{\calT}{{\mathcal T}}

\newcommand{\nbeq}{\begin{equation}}
\newcommand{\neeq}{\end{equation}}
\newcommand{\beq}{\begin{equation*}}
\newcommand{\eeq}{\end{equation*}}

\newcommand{\fp}{\mathrm{FP}}
\newcommand{\fpn}{\mathrm{FP}_{n}}
\newcommand{\fh}{\mathrm{FH}}
\newcommand{\fhn}{\mathrm{FH}_{n}}

\DeclareMathOperator{\St}{St}

\DeclareMathOperator{\BB}{BB}

\DeclareMathOperator{\fin}{fin}
\newcommand{\vfin}{V_{\fin}}
\newcommand{\vinfty}{V_{\infty}}
\newcommand{\vinftyminus}{V_{\infty}^-}

\DeclareMathOperator{\lk}{lk}
\DeclareMathOperator{\ulk}{ulk}
\DeclareMathOperator{\dlk}{dlk}

\newcommand{\aslk}{\lk^{\uparrow}}
\newcommand{\deslk}{\lk^{\downarrow}}

\newcommand{\asulk}{\ulk^{\uparrow}}
\newcommand{\desulk}{\ulk^{\downarrow}}

\newcommand{\asdlk}{\dlk^{\uparrow}}
\newcommand{\desdlk}{\dlk^{\downarrow}}

\newcommand{\flag}[1]{{#1}^{\blacktriangle}}

\newcommand\mwedge{\scalebox{.95}[1]{$\bigwedge$}}
\newcommand\mvee{\scalebox{.95}[1]{$\bigvee$}}

\DeclareMathOperator{\supp}{supp}

\begin{document}

\title[Bestvina--Brady kernels with torsion]{Finiteness properties of Bestvina--Brady kernels for graph products with torsion}

\author{Eduardo Flores Martínez}
\address{School of Mathematical Sciences, University of Southampton, Highfield, Southampton, SO17 1BJ, United Kingdom}
\email{E.Flores-Martinez@soton.ac.uk}

\author{Luis Jorge Sánchez Saldaña}
\address{Departamento de Matem\'aticas, Facultad de Ciencias, Universidad Nacional Aut\'onoma de M\'exico, 04510 Ciudad de M\'exico, Mexico}
\email{luisjorge@ciencias.unam.mx}


\begin{abstract}
We study Bestvina--Brady kernels associated to graph products whose vertex
groups are either infinite cyclic or finite. Using a $\mathrm{CAT}(0)$-cubical  complex constructed by Ruane and Witzel on which the graph product, with defining graph $\Gamma$, acts
and an equivariant Morse function, we give sufficient conditions for homological finiteness properties in
terms of the topology of the flag complex $\flag{\Gamma}$, the subcomplex
$\flag{\Gamma}_\infty$ spanned by the infinite vertices, and the ascending and descending links of certain vertices. We also prove a partial converse, showing that, under additional vanishing hypotheses on the homology in the next degree, non-vanishing homology of
these complexes obstructs higher finiteness properties of the
Bestvina--Brady kernel. The main new
features are the presence of finite stabilizers and the resulting additional
link conditions, which are expressed inside the flag complex of $\Gamma$.
When the finite vertex groups are abelian, we characterize when the kernel has
finitely many conjugacy classes of finite subgroups, describe the normalizers
of its finite subgroups as Bestvina--Brady kernels of smaller graph products,
and obtain a criterion for the kernel to be of Bredon type
$\underline{\fp}_n$.
\end{abstract}

\maketitle

\section{Introduction}\label{section:newer-introduction}

In their celebrated paper \cite{BB}, Bestvina and Brady associated to every
finite flag complex $L$ the kernel $\BB_L$ of the homomorphism
$A_L\to\dbZ$ from the right-angled Artin group $A_L$ sending each standard
generator to $1$. They showed that the finiteness properties of $\BB_L$ are
governed by the connectivity of $L$, and used these groups to produce quite interesting
examples of groups with prescribed finiteness properties. Since then
Bestvina--Brady groups have been extensively studied.

The goal of this article is to develop an analogous approach for certain graph
products. We start with a finite graph $\Gamma$ in which every vertex is
labelled with a group, called the vertex group. The corresponding graph product
is obtained from the free product of the vertex
groups by imposing commutation relations between the groups at adjacent
vertices; we recall the definition in \cref{section:graph products}.
Throughout the article, unless stated otherwise, the vertex groups are either
infinite cyclic or non-trivial finite, and there is at least one infinite
vertex (\cref{conv:standing:graph-products}). For such a graph product
$G(\Gamma)$, sending each chosen generator of an infinite cyclic vertex group
to $1\in\dbZ$ and every element of a finite vertex group to $0$ defines a
natural epimorphism $f\colon G(\Gamma)\to\dbZ$, and we denote its kernel by
$\BB_\Gamma$. We concentrate only on the homological finiteness properties
$\fp_n(R)$ of $\BB_\Gamma$. Our main results relate these properties to the
degree of connectivity of certain subcomplexes of the flag complex
$\flag{\Gamma}$ of the underlying graph $\Gamma$. We emphasize that all the
complexes appearing in our statements are subcomplexes of $\flag{\Gamma}$:
namely $\flag{\Gamma}$ itself, the full subcomplex $\flag{\Gamma}_\infty$
spanned by the infinite vertices, and, for a non-empty clique $C$ of finite
vertices, the link $\lk(C,\flag{\Gamma})$ and its full subcomplex
$\lk(C,\flag{\Gamma})_\infty$ spanned by the infinite vertices. 

We recall the theorem of Bestvina and Brady in the form that is most relevant
to this paper.

\begin{theorem}[Bestvina--Brady \cite{BB}]\label{thm:newer:classical:bb}
Let $L$ be a finite flag complex, let $R$ be an associative ring with $1$, and
let $\BB_L=\ker(A_L\to\dbZ)$. Then:
\begin{enumerate}
    \item $\BB_L$ is of type $\fp_{n+1}(R)$ if and only if $L$ is homologically
    $n$-connected;
    \item $\BB_L$ is of type $\fp(R)$ if and only if $L$ is acyclic;
    \item $\BB_L$ is finitely presented if and only if $L$ is simply connected.
\end{enumerate}
\end{theorem}

Our main results can be understood as variations of items (1) and (3) of this theorem for
$\BB_\Gamma$. The first one is a positive criterion, proved in
\cref{section:proof-main-a} as \cref{thm:principal:connectivity implies finiteness}.

\begin{theorem}\label{thm:newer:A}
Let $\Gamma$ be a labeled graph whose vertex groups are finite or infinite
cyclic and let $R$ be an associative ring with 1. Then:
\begin{enumerate}
    \item if $\flag{\Gamma}$ and $\flag{\Gamma}_\infty$ are homologically
    $n$-connected, and for every non-empty clique $C$ of finite vertices the
    complexes $\lk(C,\flag{\Gamma})$ and $\lk(C,\flag{\Gamma})_\infty$ are
    homologically $n$-connected and homologically $(n-|C|)$-connected,
    respectively, then $\BB_\Gamma$ is of type $\fp_{n+1}(R)$;
    \item if $\flag{\Gamma}$ and $\flag{\Gamma}_\infty$ are simply connected,
    if $\lk(C,\flag{\Gamma})$ is simply connected for every non-empty clique
    $C$ of finite vertices, and if $\lk(w,\flag{\Gamma})_\infty$ is non-empty and connected
    for every finite vertex $w$ while $\lk(C,\flag{\Gamma})_\infty$ is
    non-empty for every clique $C$ of finite vertices with $|C|=2$, then
    $\BB_\Gamma$ is finitely presented.
\end{enumerate}
\end{theorem}

Compared with the classical theorem, the single condition on $L$ is replaced
by a finite family of conditions on subcomplexes of $\flag{\Gamma}$, which
reflect the presence of finite vertex groups. We do not obtain type $\fp(R)$
when these complexes are acyclic, as in the second item of
\cref{thm:newer:classical:bb}, but only type $\fp_\infty(R)$; this is to be
expected, since $\BB_\Gamma$ has torsion.

To state the second result we need two conditions on a simplicial complex $L$
(here $L$ will always be one of the subcomplexes of $\flag{\Gamma}$ above). We
use augmented reduced homology, so that $\widetilde H_{-1}(\emptyset;R)=R$,
and for $q\in\dbZ$ we say that $L$ satisfies the
\emph{$(R,q)$-vanishing condition} if
\[
    \widetilde H_k(L;R)=0
    \qquad\text{for }k=q+2\text{ and for }-1\le k\le q.
\]
In other words, $L$ is homologically $q$-connected over $R$, and in addition
$\widetilde H_{q+2}(L;R)$ vanishes. We say that $L$ is
\emph{strictly $(R,q)$-acyclic} if it satisfies the $(R,q)$-vanishing
condition and, moreover, $\widetilde H_{q+1}(L;R)\neq0$. Thus a strictly
$(R,q)$-acyclic complex is one whose first non-trivial reduced homology group
appears exactly in degree $q+1$, and the next group vanishes. The precise
definitions are given in \cref{section:homology-level-sets}. 

The second result is a partial converse of the homological part of
\cref{thm:newer:A}, proved in \cref{section:homology-level-sets} as
\cref{thm:partial:converse:principal:links}.

\begin{theorem}\label{thm:newer:B}
Let $R$ be a nonzero commutative ring and let $n\geq0$. Suppose that
$\flag{\Gamma}$ and $\lk(C,\flag{\Gamma})$, for every non-empty clique $C$ of
finite vertices, satisfy the $(R,n-1)$-vanishing condition, and that
$\flag{\Gamma}_\infty$ and $\lk(C,\flag{\Gamma})_\infty$ satisfy the
$(R,n-1)$- and $(R,n-|C|-1)$-vanishing conditions, respectively. If, in
addition, at least one of the following complexes is strictly acyclic in the
indicated degree:
\begin{enumerate}
    \item $\flag{\Gamma}$ is strictly $(R,n-1)$-acyclic;
    \item $\lk(C_0,\flag{\Gamma})$ is strictly $(R,n-1)$-acyclic for some
    non-empty clique $C_0$ of finite vertices;
    \item $\flag{\Gamma}_\infty$ is strictly $(R,n-1)$-acyclic;
    \item $\lk(C_0,\flag{\Gamma})_\infty$ is strictly $(R,n-|C_0|-1)$-acyclic
    for some non-empty clique $C_0$ of finite vertices;
\end{enumerate}
then $\BB_\Gamma$ is not of type $\fp_{n+1}(R)$.
\end{theorem}

Theorems~\ref{thm:newer:A} and~\ref{thm:newer:B} do not completely
characterize the finiteness properties of $\BB_\Gamma$, but they do relate
them to the connectivity of the subcomplexes of $\flag{\Gamma}$ listed above.
The converse of \cref{thm:newer:A} is false: for the graphs
$\Gamma_k$ of \cref{section:new-examples} with $k\geq2$, the link of the maximal clique of
finite vertices is empty, so the hypotheses of \cref{thm:newer:A} fail for
every $n\geq0$, yet $\BB_{\Gamma_k}$ is of type $\fp_{k-1}(R)$
(\cref{prop:products-virtually-free}).
The role of the living subcomplex in the work of Meier, Meinert and VanWyk
\cite{MMV98} on the $\Sigma$-invariants of right-angled Artin groups, and in
the geometric proof of Bux and Gonzalez \cite{BG99}, is played here by
$\flag{\Gamma}_\infty$, while cliques of finite vertices contribute the
additional link conditions. Related results for Artin groups satisfying the
$K(\pi,1)$-conjecture appear in the work of Escartín-Ferrer and
Martínez-Pérez \cite{MPEF24}.

Our strategy is similar to that of Bestvina and Brady. The role of the
universal cover of the Salvetti complex is played by the $\mathrm{CAT}(0)$ cube complex
$X$ constructed by Ruane and Witzel, on which $G(\Gamma)$ acts properly and
cocompactly with finite stabilizers
(\cref{section:ruane-witzel-complex}); in fact $X$ is a model for
$\underline{E}G(\Gamma)$, although we do not use this. We define an
equivariant Morse function on $X$ and describe its links, ascending links and
descending links at every vertex
(\cref{section:morse-function-rw,section:ascending-descending-links}); each
of them is described in terms of the subcomplexes of $\flag{\Gamma}$
mentioned above. In contrast with the classical case, not all ascending and
descending links are isomorphic. \cref{thm:newer:A} follows directly from
the standard machinery of discrete Morse theory, after passing to
finite-index torsion-free subgroups to account for the finite stabilizers. To
prove \cref{thm:newer:B} we compute the homology of the level sets of the
Morse function. Although this is inspired by Bestvina and Brady, the argument
does not carry over verbatim, and we must impose the additional vanishing
hypotheses under which the homology of the level sets can be computed.

As an application, we study the Bredon finiteness properties of
$\BB_\Gamma$ when the finite vertex groups are abelian
(\cref{section:bredon-finiteness}). By a standard criterion, these reduce to
two conditions: having finitely many conjugacy classes of finite subgroups,
which we characterize in \cref{section:bredon-conjugacy-classes}, and having
normalizers of finite subgroups of type $\fp_n(R)$, which we study in
\cref{section:bredon-normalizers} by describing these normalizers as graph
products. We close with examples illustrating the
criteria and constructing kernels with prescribed finiteness behavior
(\cref{section:example,section:new-examples}).

The article is organized as follows. In \cref{section:preliminaries} we
recall the terminology on graph products and simplicial complexes, the
Bestvina--Brady homomorphism, and our conventions on finiteness and
connectivity. The Ruane--Witzel complex is constructed in
\cref{section:ruane-witzel-complex}, its Morse function is defined in
\cref{section:morse-function-rw}, and the links are described in
\cref{section:ascending-descending-links}. \cref{thm:newer:A} is proved in
\cref{section:proof-main-a} and \cref{thm:newer:B} in
\cref{section:homology-level-sets}. The Bredon finiteness properties and the
normalizers are treated in \cref{section:bredon-finiteness}, and the
examples are in \cref{section:example}.

We do not know to what extent the additional hypotheses of \cref{thm:newer:B}
can be removed; these hypotheses concern the homology in the degree
immediately above the one in which the connectivity of the links fails. Nor do
we know whether the converse of part~(2) of \cref{thm:newer:A} holds.

\section{Preliminaries}\label{section:preliminaries}

\subsection{Simplicial complexes and graphs}\label{section:simplicial complexes}

An \emph{abstract simplicial complex} $\mathcal{C}$ is a collection of finite subsets of a countable set of vertices $V$ such that every subset of an element of $\mathcal{C}$ is again an element of $\mathcal{C}$. An element of $\mathcal{C}$ with $n+1$ vertices is called an \emph{$n$-simplex}; in particular, vertices are the $0$-simplices. The \emph{dimension} of $\mathcal{C}$ is the supremum of the dimensions of its simplices. A \emph{map} $\mathcal{C}\to\mathcal{C}'$ of abstract simplicial complexes is a map between their vertex sets that sends simplices to simplices.

The geometric realization of an abstract simplicial complex $\calC$ is denoted by $|\calC|$. We use the term \emph{simplicial complex} for either an abstract simplicial complex or its geometric realization, identifying the two when no confusion can arise.

Let $\mathcal{C}_1,\dots, \mathcal{C}_m$ be abstract simplicial complexes. The \emph{join} $\mathcal{C}_1*\cdots *\mathcal{C}_m$ is the abstract simplicial complex whose set of vertices is the disjoint union of the vertex sets of the complexes $\mathcal{C}_i$, and whose simplices are of the form $\sigma_1\sqcup\cdots\sqcup\sigma_m$, where each $\sigma_i$ is a simplex of $\mathcal{C}_i$. We have that $|\mathcal{C}_1*\cdots *\mathcal{C}_m|$ is homeomorphic to $|\mathcal{C}_1|*\cdots *|\mathcal{C}_m|$.

A \emph{simplicial graph} is a simplicial complex of dimension at most one. A subgraph $\Gamma'\subseteq\Gamma$ is \emph{full} (or \emph{induced}) if whenever two vertices of $\Gamma'$ span an edge in $\Gamma$, they also span an edge in $\Gamma'$.

\begin{definition}\label[definition]{def:flag}
Let $\Gamma$ be a simplicial graph.
\begin{enumerate}
    \item The \emph{flag complex} associated to $\Gamma$, denoted by $\flag{\Gamma}$, is the simplicial complex obtained from $\Gamma$ by adjoining a simplex for every finite set of pairwise adjacent vertices.
    \item A \emph{clique} of $\Gamma$ is a full subgraph in which every pair of distinct vertices spans an edge.
    \item If $\mathcal{C}$ is a simplicial complex and $S$ is a simplex of $\mathcal{C}$, the \emph{link} of $S$ in $\mathcal{C}$ is
    \[
        \lk(S,\mathcal{C})=\{T\in\mathcal{C}\mid S\cap T=\emptyset\text{ and }S\cup T\in\mathcal{C}\}.
    \]
\end{enumerate}
\end{definition}

If $C$ is a clique of a simplicial graph $\Gamma$, we identify $C$ with its
set of vertices, which is a simplex of $\flag{\Gamma}$, and denote its link
in $\flag{\Gamma}$ by $\lk(C,\flag{\Gamma})$.
We denote by $\St_\Gamma(C)$ the full subgraph of $\Gamma$ spanned by the vertices of $C$ together with the vertices of $\lk(C,\flag{\Gamma})$.

Let $(X,\le)$ be a poset. The \emph{order complex} associated to $X$ is the abstract simplicial complex whose vertices are the elements of $X$ and whose simplices are the finite chains $x_0<\cdots<x_n$. The geometric realization of this order complex is denoted by $|X|$. Throughout, we write $\simeq$ both for homotopic maps and for homotopy equivalent spaces.

The following standard lemma will be used in the proof of \cref{thm:description:as:des:links}; we include it to make the relevant arguments self-contained and easier to read.

\begin{lemma}\label{lemma:poset:homotopy}
Let $f,g\colon P\to Q$ be order-preserving maps. If $f(x)\leq g(x)$ for every $x\in P$, or if the reverse inequalities hold, then the induced maps $|f|,|g|\colon |P|\to |Q|$ are homotopic. Consequently, if $i\colon Q\to P$ and $r\colon P\to Q$ are order-preserving maps such that $ri=\mathrm{id}_Q$ and $ir$ is pointwise comparable with $\mathrm{id}_P$, meaning that either $ir(x)\leq x$ for every $x\in P$ or the reverse inequalities hold, then $|i|$ and $|r|$ are homotopy inverses.
\end{lemma}
\begin{proof}
Regarding posets as categories, the pointwise inequalities define a natural transformation between $f$ and $g$, which induces a homotopy between their maps on order complexes. The final assertion follows from $|r|\circ|i|=\mathrm{id}_{|Q|}$ and $|i|\circ|r|\simeq\mathrm{id}_{|P|}$.
\end{proof}

Let $\calC$ be an abstract simplicial complex, and let $\calC^{\times}$ be the poset of non-empty simplices of $\calC$, ordered by inclusion. The \emph{barycentric subdivision} $\calC^b$ is the order complex associated to $\calC^{\times}$.

\subsection{Graph products}\label{section:graph products}

A \emph{labeled graph} consists of a simplicial graph $\Gamma$ with finite
vertex set $V$ and a \emph{vertex group} $G_v$ for every $v\in V$. The
\emph{graph product} $G(\Gamma)$ associated to a labeled graph $\Gamma$ is
the quotient of the free product of the vertex groups by the relations
$[g,h]=1$ for every $g\in G_v$ and $h\in G_w$ such that $\{v,w\}$ is an
edge of $\Gamma$.

\begin{convention}\label[convention]{conv:standing:graph-products}
Throughout the article, the vertex groups are either infinite cyclic or
non-trivial finite groups, and we assume that there is at least one infinite
vertex. 
\end{convention}

We denote by $\vfin$ and $\vinfty$ the sets of vertices labeled by
finite groups and infinite cyclic groups, respectively, so that
$V=\vfin\sqcup\vinfty$ and
$\vinfty\neq\emptyset$. The elements of these sets are called
\emph{finite vertices} and \emph{infinite vertices}, respectively. For each
$v\in\vinfty$, we fix a generator of $G_v$ and identify this generator with
$v$.

A \emph{syllable word} is a product $g_1\cdots g_n$, where each
$g_i\in G_{v_i}\setminus\{1\}$ for some $v_i\in V$. Such a word is
\emph{reduced} if no two syllables belonging to the same vertex group can be
brought together by repeatedly interchanging adjacent syllables from
commuting vertex groups. We recall the normal form theorem for graph
products.

\begin{theorem}[{\cite[Theorem~3.9]{Green90}}]
\label{thm:normal:form:graph:products}
Every element of $G(\Gamma)$ is represented by a reduced syllable word, and
any two reduced syllable words representing the same element differ by a
finite sequence of interchanges of adjacent syllables from commuting vertex
groups. In particular, the syllable length and the set of vertices occurring
in a reduced word are well defined.
\end{theorem}

If $\Lambda$ is a full labeled subgraph of $\Gamma$, we use
$\Lambda_\infty$ and $\flag{\Lambda}_\infty$ for the corresponding full
subgraph and flag subcomplex spanned by the infinite vertices of $\Lambda$.
For a clique $C$ of such a subgraph $\Lambda$, we also write
\(
    \lk(C,\flag{\Lambda})_\infty
    :=\lk(C,\flag{\Lambda})\cap\flag{\Lambda}_\infty.
\)

If $\Lambda$ is a full labeled subgraph of $\Gamma$, the graph product
$G(\Lambda)$ embeds canonically in $G(\Gamma)$; we identify $G(\Lambda)$
with its image and call it the \emph{special subgroup} associated to
$\Lambda$.

The chosen generators determine an epimorphism
\[
    f\colon G(\Gamma)\longrightarrow\dbZ
\]
that sends every infinite vertex to $1$ and every element of a finite vertex
group to $0$. We call
\[
    \BB_\Gamma:=\ker f
\]
the \emph{Bestvina--Brady kernel} associated to $\Gamma$. This homomorphism
generalizes the standard Bestvina--Brady homomorphism for right-angled Artin
groups.

For a full labeled subgraph $\Lambda$ of $\Gamma$ we put
$\BB_\Lambda:=\ker(f|_{G(\Lambda)})=G(\Lambda)\cap\BB_\Gamma$. This makes
sense even when $\Lambda$ has no infinite vertices, so that $\Lambda$ does
not satisfy \cref{conv:standing:graph-products}; in that case
$f|_{G(\Lambda)}$ is trivial and $\BB_\Lambda=G(\Lambda)$. This applies, for
instance, to $\BB_{\St_\Gamma(\calC)}$ when
$\lk(\calC,\flag{\Gamma})_\infty=\emptyset$.

\subsection{Bestvina--Brady Morse theory}\label{section:bb-morse-theory}

A \emph{Morse function} on a simplicial complex $X$ is a continuous map $h\colon X\to \dbR$ such that $h(X^{(0)})$ is a discrete subset of $\dbR$ and whose restriction to each simplex is affine and attains its maximum and minimum at unique vertices.

\begin{definition}\label[definition]{def:asclk}
Let $v$ be a vertex of $X$. The \emph{ascending link} $\aslk v$ with respect to a Morse function $h$ on $X$ is the subcomplex of $\lk(v,X)$ consisting of those simplices $S$ such that $h$ attains its minimum on the simplex $S*\{v\}$ uniquely at $v$. The \emph{descending link} $\deslk v$ is defined analogously, requiring that $h$ attain its maximum on $S*\{v\}$ uniquely at $v$.
\end{definition}

For $J\subseteq\mathbb{R}$, write $X_J=h^{-1}(J)$. The topology of $X$ is described by a Morse function and its associated descending links as shown by the \emph{Morse lemma} below.

\begin{lemma}[{\cite{BB}}]
    Let $h:X\to\mathbb{R}$ be a Morse function. Suppose that $J\subset J'\subset\mathbb{R}$ are closed and connected, $\inf J=\inf J'$, and $J'\setminus J$ contains exactly one value $r\in h(X^{(0)})$. Then $X_{J'}$ is homotopy equivalent to the complex obtained from $X_J$ by coning off the descending link $\deslk v$ of every vertex $v$ with $h(v)=r$.
\end{lemma}

There is an analogous statement for ascending links when $\sup J=\sup J'$.

\subsection{Connectivity conventions}

Let $R$ be a nonzero associative ring with unit. We use the convention that
$\widetilde H_{-1}(X;R)=0$ if and only if $X$ is non-empty, and that
$\widetilde H_i(X;R)=0$ for every $i<-1$. A space $X$ is \emph{homologically
$n$-connected over $R$} if $\widetilde H_i(X;R)=0$ for all $i\le n$. Thus,
homological $(-1)$-connectivity means nonemptiness, while homological
$n$-connectivity imposes no condition when $n\leq -2$. The space $X$ is
\emph{$R$-acyclic} if $\widetilde H_i(X;R)=0$ for every $i\geq -1$. When the
coefficient ring is clear from the context, we simply say that $X$ is
homologically $n$-connected (respectively, acyclic).

\subsection{Finiteness properties of groups}

We recall the homological finiteness properties used throughout the article.

\begin{definition}
Let $R$ be a nonzero associative ring with unit, and let $G$ be a group.
\begin{enumerate}
    \item The group $G$ is of type $\fhn(R)$ if it acts freely, properly discontinuously, cellularly, and cocompactly on a CW-complex $X$ such that $\widetilde{H}_k(X;R)=0$ for every $k<n$. If such an action exists on an $R$-acyclic CW-complex, then $G$ is of type $\fh(R)$.
    \item The group $G$ is of type $\fpn(R)$ if the trivial $RG$-module $R$ admits a projective resolution whose modules are finitely generated through dimension $n$. It is of type $\fp(R)$ if $R$ admits a finite-length resolution by finitely generated projective $RG$-modules.
\end{enumerate}
\end{definition}

It is standard that a group of type $\fhn(R)$ is of type $\fpn(R)$, and that type $\fh(R)$ implies type $\fp(R)$.

\section{The Ruane--Witzel complex}

Throughout this section, we fix a labeled graph $\Gamma$ satisfying \cref{conv:standing:graph-products}. Ruane and Witzel \cite{RW16} constructed a simplicial complex $X$ on which $G(\Gamma)$ acts simplicially, properly, and cocompactly. We now recall their construction.

\subsection{The doubled graph and brackets of cliques}\label{section:brackets of cliques}

To define the Ruane--Witzel complex associated to a graph product we construct the following simplicial graph from $\Gamma$ by doubling infinite vertices.

\begin{definition}\label[definition]{def:delta}
     Let $V$ be the set of vertices of $\Gamma$.
\begin{itemize}
    \item Let $\vinftyminus$ be a set of formal inverses of infinite vertices, disjoint from $V$, and write $\vinftyminus=\{v^{-1}\colon v\in \vinfty \}$. We identify each formal symbol $v^{-1}$ with the inverse of the chosen generator $v\in G_v$ and set $(v^{-1})^{-1}=v$.
    \item  Denote by $\hat V$ the disjoint union $\vfin\sqcup\vinfty\sqcup \vinftyminus$. Let $\hat V\to V$ be the function that is the identity on $\vfin$, and the natural 2-to-1 projection $\vinfty\sqcup \vinftyminus\to \vinfty$.
    \item Define a new simplicial graph $\Delta$ whose vertex set is $\hat V$ and such that two vertices are adjacent if their images under the projection $\hat V\to V$ are adjacent in $\Gamma$. In particular, for $v\in \vinfty$, $v$ and $v^{-1}$ are not adjacent in $\Delta$.
    \item The elements of $\vinfty$ (thought of as vertices of $\Delta$) are called \emph{positive vertices}, and the elements of $\vinftyminus$ are called \emph{negative vertices}.
\end{itemize}

\end{definition}

We use the vertices of $\Delta$ to define the basic building blocks of the Ruane--Witzel complex.

\begin{definition}
Let $v\in \hat V$. The \emph{bracket} $[v]$ is defined as
\[
[v]=\begin{cases}
G_v & \text{if }v\in \vfin,\\
\{1,v\} & \text{otherwise.}
\end{cases}
\]
We identify each clique with its vertex set. The bracket of a clique $C=\{v_1,\dots,v_k\}$ is $[C]=[v_1]\cdots[v_k]$, where the product is the setwise product in $G(\Gamma)$. By convention, $[\emptyset]=\{1\}$.
\end{definition}

For $g\in G(\Gamma)$ and a clique $C$ of $\Delta$, we write
\[
g[C]=\{gx\mid x\in [C]\}
\]
for the left $g$-translate of the bracket $[C]$.

The ordering of $v_1,\dots,v_k$ is irrelevant because the corresponding vertex groups commute. Moreover, $[C]$ is always a finite subset of $G(\Gamma)$. By \cref{thm:normal:form:graph:products}, every element of $[C]$ has a unique expression
\[
g_1\cdots g_k,
\qquad g_i\in[v_i].
\]
Consequently, the multiplication map
\[
\prod_{i=1}^k[v_i]\longrightarrow[C],
\qquad (g_1,\dots,g_k)\longmapsto g_1\cdots g_k,
\]
is a bijection. For $g\in[C]$ and $v\in C$, we call the corresponding entry of this tuple the \emph{$v$-coordinate} of $g$ and denote it by $g_v$.

Let $C_{\fin}=C\cap\vfin$ and $C_\infty=C\cap(\vinfty\sqcup\vinftyminus)$. The bracket $[C]$ factors as $[C]=[C_{\fin}][C_\infty]$. Thus every $g\in[C]$ has a unique factorization $g=g_{\fin}g_\infty$, where $g_{\fin}\in[C_{\fin}]$ and $g_\infty\in[C_\infty]$. If $C_\infty=\{v_1,\dots,v_l\}$, then
\[
g_\infty=v_1^{\epsilon_1}\cdots v_l^{\epsilon_l}
\qquad\text{for uniquely determined }\epsilon_i\in\{0,1\}.
\]

\begin{definition}
Let $C$ be a clique in $\Delta$ and let $W$ be a subset of $C_\infty$. The \emph{$W$-reflection of $C$}, denoted $C\vert W$, is the clique with vertex set
\[
C\vert W=(C\setminus W)\cup\{v^{-1}\colon v\in W\}.
\]
This is a clique because $v$ and $v^{-1}$ have the same image in $\Gamma$.
Let $g\in[C]$. The \emph{$g$-reflection of $C$} is the $\{v\in C_\infty\colon g_v\neq 1\}$-reflection of $C$, and it is denoted by $C\vert g$.
\end{definition}

We end this subsection with the following basic properties of brackets of cliques.

\begin{lemma}\label[lemma]{lemma:properties:brackets}
Let $C$ and $D$ be cliques in
$\Delta$, and let $g\in G(\Gamma)$. Then:
\begin{enumerate}
    \item\label{lemma:properties:brackets:subset} $[C]\subseteq[D]$ if and only if $C\subseteq D$.
    \item\label{lemma:properties:brackets:equality} $[C]=[D]$ if and only if $C=D$.
    \item $[C_{\fin}]$ is a finite subgroup of $G(\Gamma)$.
    \item $1\in g[C]$ if and only if $g^{-1}\in[C]$. In this case,
    \(
        g[C]=[C\vert g^{-1}].
    \)
    \item $g[C]=[C]$ if and only if $g\in[C_{\fin}]$.
\end{enumerate}
\end{lemma}
\begin{proof}
If $C\subseteq D$, then $[C]\subseteq[D]$. Conversely, if
$[C]\subseteq[D]$, \cref{thm:normal:form:graph:products,conv:standing:graph-products}
show that every vertex of $C$, with its orientation when it is
infinite, belongs to $D$. This proves (1), and (2) follows immediately.

Since the finite vertices of $C$ form a clique,
$[C_{\fin}]=\prod_{v\in C_{\fin}}G_v$ is a finite subgroup, proving (3).

Clearly, $1\in g[C]$ if and only if $g^{-1}\in[C]$. In this case,
translation by $g$ permutes each finite coordinate and replaces
$\{1,v\}$ by $\{1,v^{-1}\}$ precisely when $(g^{-1})_v\neq1$. Hence
$g[C]=[C\vert g^{-1}]$, proving (4).

If $g\in[C_{\fin}]$, then $g[C]=[C]$ because $[C_{\fin}]$ is a subgroup
commuting with $[C_\infty]$. Conversely, if $g[C]=[C]$, then (4) and (2)
give $C=C\vert g^{-1}$. Thus $g^{-1}$ has no non-trivial infinite coordinates, so
$g^{-1}\in[C_{\fin}]$, and hence $g\in[C_{\fin}]$. This proves (5).
\end{proof}

\subsection{The Ruane--Witzel complex}\label{section:ruane-witzel-complex}

Let $\calS$ be the poset of brackets $[C]$, where $C$ ranges over the cliques
of $\Delta$, ordered by inclusion as subsets of $G(\Gamma)$.
Define $K$ to be the geometric realization of $\calS$, that is, $K$ is a simplicial complex whose vertices are labeled by brackets $[C]$ and an $n$-simplex is given by a \emph{flag} $[C_0]\subsetneq \cdots \subsetneq [C_n]$.

Let $\calT=\{g[C]\mid g\in G(\Gamma),\ C\text{ is a clique of }\Delta\}$ be the poset of translates of elements of $\calS$, again ordered by inclusion. Define $X$ to be the geometric realization of $\calT$. Thus the vertices of $X$ are the translates $g[C]$, with $g\in G(\Gamma)$ and $C$ a clique of $\Delta$, and its simplices are the finite chains
\(
    g_0[C_0]\subsetneq g_1[C_1]\subsetneq \cdots \subsetneq g_n[C_n]
\)
in this poset. The group $G(\Gamma)$ acts on $\calT$ by left multiplication, $h\cdot g[C]=hg[C]$, and therefore acts simplicially on $X$.

\begin{lemma}[{\cite{RW16}}]\label[lemma]{lemma:basic-properties-rw}
We have the following basic properties of $K$ and $X$.
\begin{enumerate}
    \item With its standard piecewise Euclidean metric, $X$ is a locally finite $\mathrm{CAT}(0)$ complex.
    \item The simplicial $G(\Gamma)$-action defined above is by isometries.
    \item The subspace $K$ is a weak fundamental domain for this action, in the sense of \cite{RW16}.
    \item The stabilizer of $[C]$ in $G(\Gamma)$ is the finite subgroup $[C_{\fin}]$. Hence every simplex stabilizer, being contained in a vertex stabilizer, is finite.
\end{enumerate}
\end{lemma}

The following lemma records useful computations in $\calT$ that will be used to describe links in $X$.

\begin{lemma}\label[lemma]{lemma:intervals:initial:endpoint:S}
\begin{enumerate}
    \item If $[C]\subseteq g[D]$, then $g^{-1}\in [D]$ and $g[D]=[D\vert g^{-1}]$.
    \item If $g[C]\subseteq [D]$, then $g\in [D]$, and there are a clique $C'\subseteq D$ and an element $a\in[D\setminus C']$ such that $g[C]=a[C']$. Moreover, $C'$ is a reflection of $C$.
    \item If $a[C']=b[E]$, where $C',E\subseteq D$, $a\in[D\setminus C']$ and $b\in[D\setminus E]$, then $C'=E$ and $a=b$.
\end{enumerate}
\end{lemma}
Thus every $g[C]\subseteq[D]$ has a unique expression $a[C']$ as in (2); we call it the \emph{normalized representative} of $g[C]$ in $[D]$.
\begin{proof}
\begin{enumerate}
    \item Since $1\in [C]$ and $[C]\subseteq g[D]$, we get $1\in g[D]$. Now the claim is an immediate consequence of \cref{lemma:properties:brackets}(4).
    
    \item Since $1\in [C]$ and $g[C]\subseteq [D]$, we have $g\in [D]$. Moreover, $[C]\subseteq g^{-1}[D]$, so (1) gives $g^{-1}[D]=[D\vert g]$, and hence $C\subseteq D\vert g$ by \cref{lemma:properties:brackets}(1). Let $W$ be the set of vertices of $C$ not in $D$. Each $u\in W$ lies in $(D\vert g)\setminus D$, so $u=v^{-1}$ for some infinite $v\in D$ with $g_v=v$. Thus $W\subseteq C_\infty$, and $C'=C\vert W$ is a reflection of $C$ contained in $D$; note that $C'=(C\cap D)\sqcup W^{-1}$, where $W^{-1}=\{u^{-1}\colon u\in W\}$.

    Let $h=\prod_{u\in W}u\in[C]$. By \cref{lemma:properties:brackets}(4), $h^{-1}[C]=[C\vert h]=[C']$, so $g[C]=gh[C']$. Since the vertex groups of $D$ commute and $g_v=v$ for $v\in W^{-1}$, we get $gh=\prod_{v\in D\setminus W^{-1}}g_v$. If $u\in C\cap D$ is infinite, then $u$ is not reflected in $D\vert g$, so $g_u=1$. Hence $gh=ac$, where $a=\prod_{v\in D\setminus C'}g_v\in[D\setminus C']$ and $c=\prod_{u\in C'_{\fin}}g_u\in[C'_{\fin}]$. By \cref{lemma:properties:brackets}(5), $c[C']=[C']$, and therefore $g[C]=a[C']$.

    \item Since $1\in[E]=b^{-1}a[C']$, \cref{lemma:properties:brackets}(4) shows that $E=C'\vert a^{-1}b$, so $E$ is a reflection of $C'$. Both $C'$ and $E$ are contained in $D$. If $C'\neq E$, then $D$ would contain an infinite vertex and its reflection, contradicting that $D$ is a clique. Thus $C'=E$. Now $b^{-1}a$ stabilizes $[C']$, so \cref{lemma:properties:brackets}(5) gives $b^{-1}a\in[C'_{\fin}]$. On the other hand, $b^{-1}a=\prod_{v\in D\setminus C'}b_v^{-1}a_v$, and the vertices of $D\setminus C'$ and of $C'$ have disjoint images in $\Gamma$. By \cref{thm:normal:form:graph:products}, $a_v=b_v$ for every $v\in D\setminus C'$, that is, $a=b$.
\end{enumerate}
\end{proof}

\section{A Morse function on the Ruane--Witzel complex}\label{section:morse-function-rw}

Recall that $\calT$ is the poset of translates $g[C]$, with $g\in G(\Gamma)$ and $C$ a clique of $\Delta$, ordered by inclusion. The Ruane--Witzel complex $X$ is the geometric realization of $\calT$. A map on the vertices $\calT$ extends uniquely to a map $X\to\dbR$ that is affine on each simplex.

Recall the homomorphism $f\colon G(\Gamma)\to\dbZ$ and the Bestvina--Brady
kernel $\BB_\Gamma=\ker f$ defined in \cref{section:graph products}.
We write
$v^+=v$ and $v^-=v^{-1}$ for the two vertices of $\Delta$ associated to an
infinite vertex $v$ of $\Gamma$; thus vertices of the form $v^+$ are positive
and vertices of the form $v^-$ are negative. Also, for a finite $v\in V$ we
have $f(x)=0$ for all $x\in [v]$.

Define
\[
F_1(g[C])=(\dim(\mathcal{S})+1)\max\{f(x)\colon x\in g[C]\},
\qquad F_2(g[C])=-|C|.
\]
The first map depends only on the subset $g[C]\subseteq G(\Gamma)$. If $g[C]=h[D]$, then \cref{lemma:properties:brackets}(4) implies that $D$ is a reflection of $C$, so $|C|=|D|$; hence $F_2$ is also well defined. Set
\[
F=F_1+F_2\colon\calT\longrightarrow\dbZ,
\]
and extend $F$ affinely over every simplex of $X$.

\begin{lemma}\label[lemma]{thm:F:fequivariant}
The function $F$ is $(\dim(\mathcal{S})+1)f$-equivariant; that is,
\[
F(gh[C])=(\dim(\mathcal{S})+1)f(g)+F(h[C])
\]
for all $g\in G(\Gamma)$ and $h[C]\in\calT$.
\end{lemma}

\begin{proof}
For $g\in G(\Gamma)$, the homomorphism property of $f$ gives
\[
F_1(gh[C])=(\dim(\mathcal{S})+1)f(g)+F_1(h[C]).
\]
Moreover, $F_2(gh[C])=-|C|=F_2(h[C])$. Adding these identities proves the result.
\end{proof}

\begin{theorem}\label{thm:F:morsesimplicial}
Let $\sigma$ be a simplex represented by a flag
$[C_0]\subsetneq\cdots\subsetneq[C_n]$. Among the $C_i$ with the largest
number of positive vertices, let $C_M$ have the smallest number of vertices.
Among those with the smallest number of positive vertices, let $C_m$ have the
largest number of vertices. Then $F$ attains its unique maximum on $\sigma$ at
$[C_M]$ and its unique minimum at $[C_m]$. In particular, $F$ is a Morse
function.
\end{theorem}
\begin{proof}
It suffices to consider a simplex in $K$, since the equivariance in
\cref{thm:F:fequivariant} preserves the vertices at which the extrema occur.
For a clique $C$, the quantity $\max\{f(x):x\in[C]\}$ is the number of
positive vertices of $C$. If $C_i$ has fewer positive vertices than $C_M$,
the factor $\dim(\mathcal{S})+1$ ensures that $F([C_i])<F([C_M])$, since
$|C|<\dim(\mathcal{S})+1$ for every clique $C$. If they have the same number
of positive vertices, the choice of $C_M$ gives $|C_M|<|C_i|$ and again
$F([C_i])<F([C_M])$. Thus $[C_M]$ is the unique maximum. The argument for
$[C_m]$ is analogous. Since $F(\calT)\subseteq\dbZ$ and the extension is
affine on every simplex, $F$ satisfies the definition of a Morse function.
\end{proof}

The following examples realize $F$ as a height function.
\begin{example}
The Ruane--Witzel complex of the connected graph with vertex groups $\mathbb{Z}$ and $\mathbb{Z}/2\mathbb{Z}$ is an infinite band. We compute the values of $F$ on the vertices of the fundamental domain (shown in blue) and its $w$-translate, where $w$ is the finite vertex and $v$ the infinite vertex.
\[\begin{matrix}
    F([v^-])=-1, & F([\emptyset])=0, & F([v^+])=2,\\
    F([v^-,w])=-2, & F([w])=-1, & F([v^+,w])=1,\\
    F(w[v^-])=-1, & F(w[\emptyset])=0, & F(w[v^+])=2.
\end{matrix}\]

Thus, we obtain the following diagram.

\begin{figure}[h]
    \centering
    \includegraphics[width=.8\linewidth]{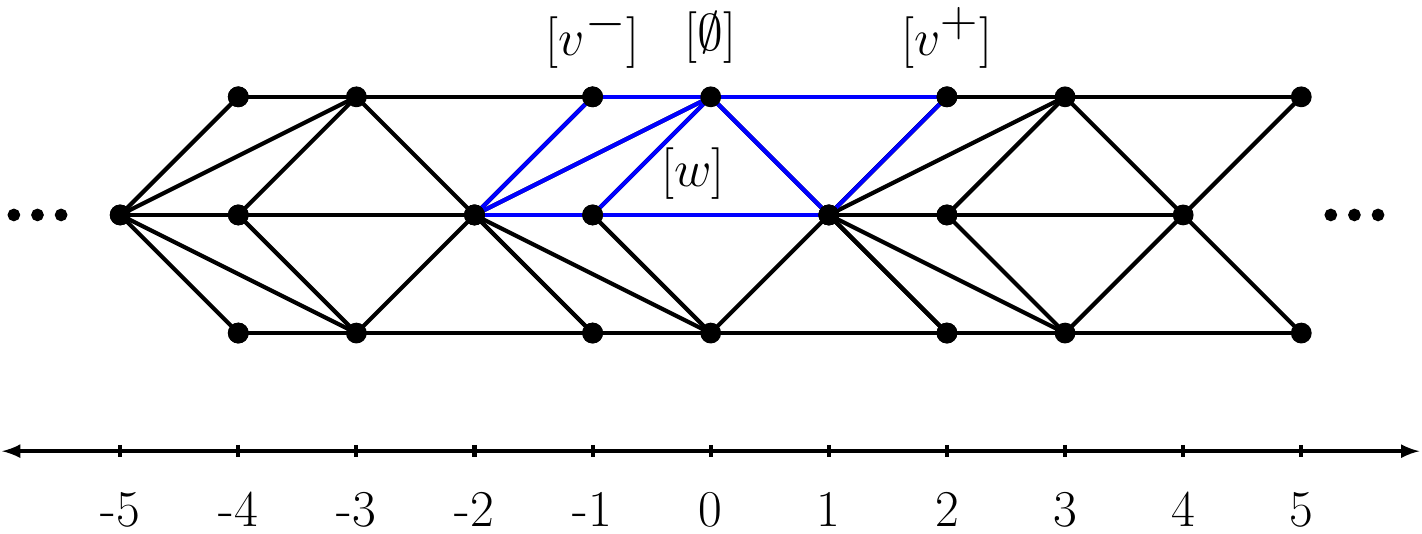}
    \caption{Example of Morse function.}
\end{figure}
\end{example}
\begin{example}
We proceed similarly for the case of the connected graph with two infinite vertices, denoted by $v$ and $w$. The following is the computation of $F$ restricted to the vertices of the fundamental domain.
\[
\begin{aligned}
F([v^+,w^+])&=4, & F([w^+])=F([v^+])&=2,\\
F([v^-,w^+])=F([v^+,w^-])&=1, & F([\emptyset])&=0,\\
F([v^-])=F([w^-])&=-1, & F([v^-,w^-])&=-2.
\end{aligned}
\]

By the equivariance of $F$ we obtain the following picture.

\begin{figure}[h]
    \centering
    \includegraphics[width=0.715\linewidth]{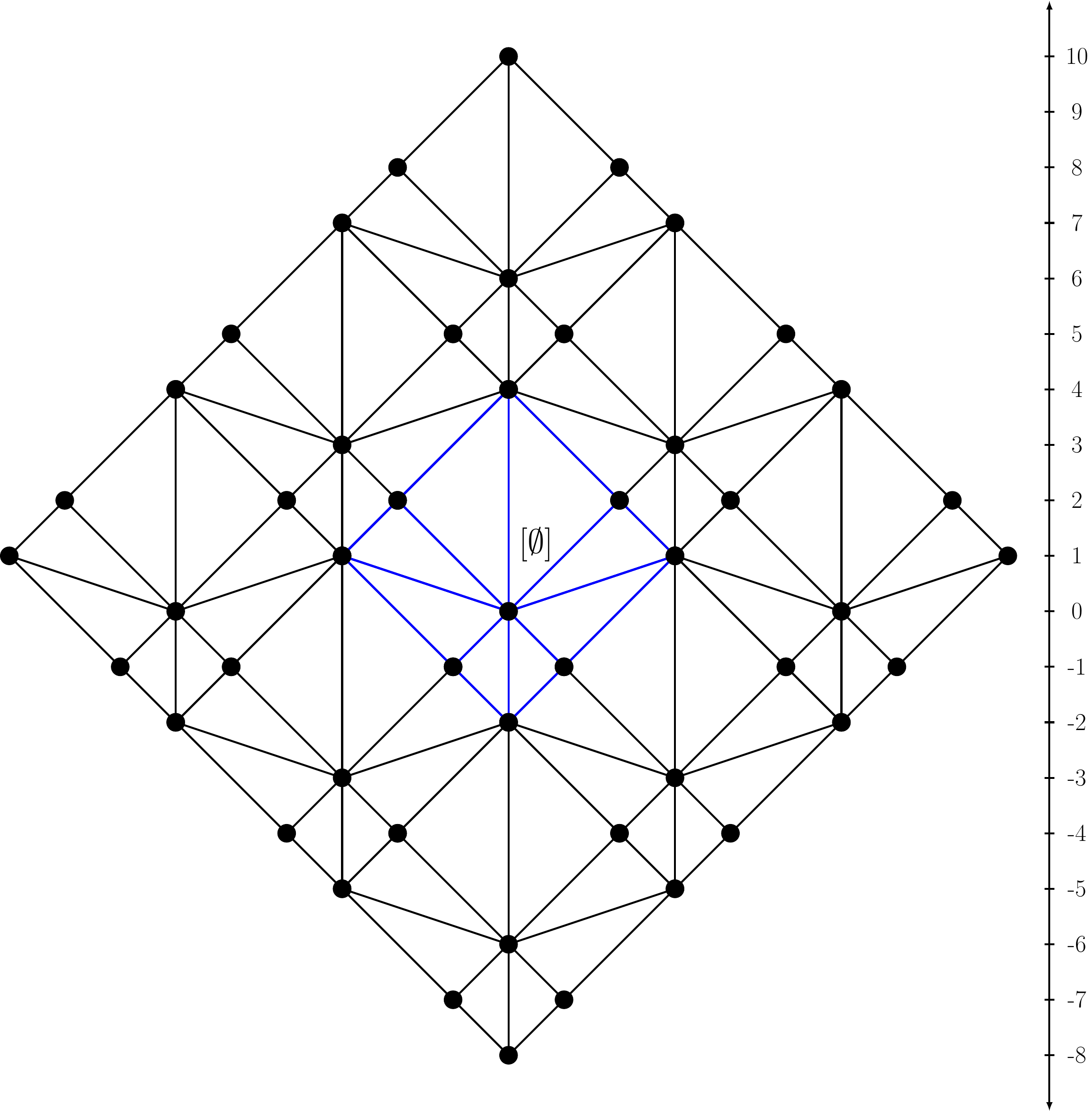}
    \caption{Morse function in the Bestvina--Brady case.}
\end{figure}

\end{example}

\section{Description of the ascending and descending links}\label{section:ascending-descending-links}

\subsection{Description of links} We first decompose the link of a bracket into its up- and down-links, and then describe the ascending and descending links associated to $F$.

Recall that every vertex of $X$ has the form $g[C]$, where $g\in G(\Gamma)$ and $C$ is a clique in $\Delta$. The action of $g^{-1}$ induces an isomorphism from the link of $g[C]$ to the link of $[C]$, so it suffices to consider vertices of the form $[C]$.

\begin{definition}\label[definition]{def:uplink}
    The \emph{up-link} $\ulk [C]$ (resp. the \emph{down-link} $\dlk [C]$) is the subcomplex of $\lk([C],X)$ spanned by the vertices $g[D]$ with $[C]\subsetneq g[D]$ (resp. $g[D]\subsetneq [C]$).
\end{definition}

The identifications below appear in \cite[Observations~3.5 and~3.6]{RW16} without proof. We include a proof for completeness; the join decomposition is the standard link decomposition for an order complex.

\begin{theorem}\label{thm:description:lk}
We have the following.
\begin{enumerate}
    \item The up-link $\ulk g[C]$ is isomorphic to  $\lk(C,\flag{\Delta})^b$.
    \item The down-link $\dlk g[C]$ is isomorphic to $(\Asterisk_{v\in C}[v])^b$, where each $[v]$ is regarded as a discrete simplicial complex.
    \item The link of a vertex of $X$ of the form $g[C]$ decomposes as
    \[
        \lk g[C]=\ulk g[C]*\dlk g[C].
    \]
\end{enumerate}
\end{theorem}
\begin{proof}
Without loss of generality we may work with $[C]$ instead of $g[C]$.

\begin{enumerate}
    \item By \cref{lemma:properties:brackets}(4), every vertex $a[D]$ of $\ulk[C]$ can be written as $[D\vert a^{-1}]$, with $C\subsetneq D\vert a^{-1}$. Define
    \begin{equation}\label{eq:iso:ulk}
    \varphi\colon \ulk[C]\to \lk(C,\flag{\Delta})^b \end{equation}
    by $\varphi(a[D])=(D\vert a^{-1})\setminus C$. Its inverse sends a non-empty clique $A$ in $\lk(C,\flag{\Delta})$ to $[C\cup A]$. These maps preserve inclusions and hence induce the required isomorphism.
    
    \item By \cref{lemma:intervals:initial:endpoint:S}(2) and~(3), every vertex of $\dlk[C]$ has a unique normalized representative $a[D]$ in $[C]$, with $D\subseteq C$ and $a\in[C\setminus D]$; moreover $D\subsetneq C$, since $D=C$ would force $a=1$. Define
    \begin{equation}\label{eq:iso:dlk}
    \psi\colon\dlk[C]\to(\Asterisk_{v\in C}[v])^b
    \end{equation}
    by $\psi(a[D])=\{(w,a_w):w\in C\setminus D\}$, where $a[D]$ is the normalized representative and the first coordinate records the copy of $[w]$ in the join. This is well defined by \cref{lemma:intervals:initial:endpoint:S}(3). Conversely, a non-empty simplex $\{(w,x_w):w\in W\}$ determines $(\prod_{w\in W}x_w)[C\setminus W]$. These assignments are inverse, and inclusions correspond to reverse inclusions of simplices, so they induce an isomorphism of order complexes.
    \item In the order complex of a poset, the link of a vertex is the join of the order complexes of its strict lower and upper intervals. These are precisely $\dlk[C]$ and $\ulk[C]$, respectively.
\end{enumerate}
\end{proof}

\subsection{Morse ascending and descending links}

We now refine the up- and down-links using the Morse function $F$ from \cref{section:morse-function-rw}. For a vertex $g[C]$, set
\[
\asulk g[C]=\aslk g[C]\cap\ulk g[C],\qquad
\asdlk g[C]=\aslk g[C]\cap\dlk g[C],
\]
and define $\desulk g[C]$ and $\desdlk g[C]$ analogously using $\deslk g[C]$. The next theorem describes these subcomplexes and hence the ascending and descending links.

\begin{definition}\label[definition]{def:subcomplexes-delta}
Let $\flag{\Delta}_{\leq0}$ (resp. $\flag{\Delta}_{\geq0}$) be the full subcomplex of $\flag{\Delta}$ spanned by the finite and negative vertices (resp. finite and positive vertices). For $\Sigma\subseteq\flag{\Delta}$, set $\Sigma_{\leq0}=\Sigma\cap\flag{\Delta}_{\leq0}$ and $\Sigma_{\geq0}=\Sigma\cap\flag{\Delta}_{\geq0}$, with analogous notation for $<0$ and $>0$.
\end{definition}

The projection $p\colon\Delta\to\Gamma$ induces simplicial isomorphisms $\flag{\Delta}_{\leq0}\cong\flag{\Gamma}$ and $\flag{\Delta}_{\geq0}\cong\flag{\Gamma}$. Since $p$ restricts to the identity on $\Delta_{\fin}$, we identify $\Delta_{\fin}$ with the full subgraph of $\Gamma$ spanned by the finite vertices. In particular, for every clique $C\subseteq\Delta_{\fin}$, we write $C=p(C)$ by abuse of notation.

\begin{observation}\label[observation]{obs:delta:gamma:finite:links}
Let $C\subseteq\Delta_{\fin}$ be a clique. Then $p$ restricts to simplicial isomorphisms
\[
    \lk(C,\flag{\Delta})\cap\flag{\Delta}_{\le0}
    \cong \lk(C,\flag{\Gamma})
    \quad\text{and}\quad
    \lk(C,\flag{\Delta})\cap\flag{\Delta}_{>0}
    \cong \lk(C,\flag{\Gamma})_\infty.
\]
\end{observation}

\begin{theorem}\label{thm:description:as:des:links}
For every clique $C$ of $\Delta$, the following hold.
\begin{enumerate}
    \item The up-links satisfy
    \[
        \desulk[C]\cong
        \lk(C,\flag{\Delta})^b\cap(\flag{\Delta}_{\leq0})^b,
    \]
    while $\asulk[C]$ is homotopy equivalent to
    $\lk(C,\flag{\Delta})^b\cap(\flag{\Delta}_{>0})^b$.
    
    \item The down-link $\desdlk[C]$ is empty if $C_\infty=\emptyset$ and is
    homotopy equivalent to $\flag{C}_\infty$ otherwise. Moreover,
    \[
        \asdlk[C]\cong
        \left(\flag{C}_\infty*\Asterisk_{v\in C_{\fin}}[v]\right)^b.
    \]
    Consequently, both down-links are contractible when $C_\infty\neq\emptyset$.
    
    \item The full links decompose as
    \[
        \deslk[C]\cong\desulk[C]*\desdlk[C],
        \qquad
        \aslk[C]\cong\asulk[C]*\asdlk[C].
    \]
\end{enumerate}
\end{theorem}
\begin{proof}
\begin{enumerate}
    \item By \cref{thm:F:morsesimplicial}, under the isomorphism $\varphi$ from \eqref{eq:iso:ulk}, the descending up-link corresponds exactly to the non-empty simplices of $\lk(C,\flag{\Delta})$ containing only finite and negative vertices. This gives the stated isomorphism.

    By \cref{thm:F:morsesimplicial}, the ascending up-link corresponds to the poset $P$ of non-empty simplices $A$ of $\lk(C,\flag{\Delta})$ with $A_{>0}\neq\emptyset$. Let $Q\subseteq P$ be the subposet of positive simplices, and let $i\colon Q\to P$ be the inclusion. The map $r\colon P\to Q$ given by $r(A)=A_{>0}$ is order-preserving and satisfies $ri=\mathrm{id}_Q$, while $ir(A)\subseteq A$. Thus \cref{lemma:poset:homotopy} gives
    \[
        |P|\simeq |Q|
        =\left(\lk(C,\flag{\Delta})\cap\flag{\Delta}_{>0}\right)^b,
    \]
    proving (1).

    \item Let $g[D]\subsetneq[C]$ be written as its normalized representative in $[C]$ (\cref{lemma:intervals:initial:endpoint:S}(2)--(3)), so that $D\subsetneq C$, $g\in[C\setminus D]$, and $g=\prod_{w\in C\setminus D}g_w$ with $g_w\in[w]$ uniquely determined. Multiplying by $g^{-1}$ and using the equivariance of $F$ shows that $g[D]$ lies in $\asdlk[C]$ exactly when $(C\vert g)\setminus D$ contains only finite and negative vertices. Equivalently, for every infinite $w\in C\setminus D$, one has $g_w=1$ when $w$ is negative and $g_w=w$ when $w$ is positive. Hence the isomorphism $\psi$ from \eqref{eq:iso:dlk} restricts to
    \[
        \asdlk[C]\cong
        \left(\flag{C}_\infty*\Asterisk_{v\in C_{\fin}}[v]\right)^b.
    \]

    If $C_\infty=\emptyset$, every vertex of $\dlk[C]$ is ascending, so $\desdlk[C]$ is empty. Otherwise, under $\psi$, the descending down-link is the order complex of the poset $P$ of simplices containing at least one opposite infinite choice: $1_v$ for $v\in C_\infty^+$ or $v$ for $v\in C_\infty^-$. Let $Q$ be the poset of non-empty faces of the simplex spanned by these opposite choices, and let $i\colon Q\to P$ be the inclusion. Intersecting with this simplex defines an order-preserving map $r\colon P\to Q$ such that $ri=\mathrm{id}_Q$ and $ir(S)\subseteq S$. By \cref{lemma:poset:homotopy}, $\desdlk[C]$ is homotopy equivalent to $|Q|$, the barycentric subdivision of $\flag{C}_\infty$. Since $\flag{C}_\infty$ is a non-empty simplex when $C_\infty\neq\emptyset$, both down-links are then contractible.

    \item By \cref{thm:description:lk}(3), the full link is the join of its up- and down-links. Taking the induced subcomplexes on vertices above, respectively below, the value $F([C])$ gives the two asserted decompositions.
\end{enumerate}
\end{proof}

\begin{corollary}\label[corollary]{cor:contractible:links}
Let $C$ be a clique of $\Delta$.
\begin{enumerate}
\item If $C=\emptyset$ then $\deslk[C]\simeq\flag{\Gamma}$ and $\aslk[C]\simeq\flag{\Gamma}_\infty$.
\item If $C_\infty\neq \emptyset$, then both $\deslk[C]$ and $\aslk[C]$ are contractible.
\item If $C_\infty=\emptyset$ and $C_{\fin}\neq\emptyset$, then $\deslk[C]\cong\lk(C,\flag{\Gamma})^b$ and $\aslk[C]$ is homotopy equivalent to the join of the barycentric subdivision of $\lk(C,\flag{\Gamma})_\infty$ with a finite wedge of spheres of dimension $|C|-1$.
\end{enumerate}
\end{corollary}
\begin{proof}
    \begin{enumerate}
        \item We have that $\desdlk[\emptyset]$ and $\asdlk[\emptyset]$ are empty and thus \begin{align*}\deslk[\emptyset]&\cong\desulk[\emptyset]\simeq\flag{\Delta}_{\leq 0}\cong\flag{\Gamma},\\
        \aslk[\emptyset]&\cong\asulk[\emptyset]\simeq\flag{\Delta}_{>0}\cong\flag{\Gamma}_\infty.
        \end{align*}
        \item If $C_\infty\neq\emptyset$, then $\flag{C}_\infty$ is a non-empty simplex. Hence both down-links are contractible by \cref{thm:description:as:des:links}(2), and the join decompositions in \cref{thm:description:as:des:links}(3) show that $\deslk[C]$ and $\aslk[C]$ are contractible.
        \item If $C_\infty=\emptyset$, then \cref{thm:description:as:des:links}(2)--(3) gives
        \[
            \deslk[C]\cong\desulk[C],
            \qquad
            \aslk[C]\cong
            \asulk[C]*\left(\Asterisk_{v\in C_{\fin}}[v]\right)^b.
        \]
        By \cref{thm:description:as:des:links}(1) and \cref{obs:delta:gamma:finite:links},
        \[
            \desulk[C]\cong\lk(C,\flag{\Gamma})^b,
            \qquad
            \asulk[C]\simeq\lk(C,\flag{\Gamma})_\infty^b.
        \]
        Finally, the join of the $|C|$ non-trivial finite sets $[v]=G_v$ is homotopy equivalent to a finite wedge of spheres of dimension $|C|-1$, which proves (3).
    \end{enumerate}
\end{proof}

\subsection{Connectivity properties of links}\label{section:link connectivity}

We now determine the connectivity and the reduced homology of the ascending and descending links associated to non-empty cliques in terms of the links of these cliques in $\flag{\Gamma}$ and in its infinite subcomplex $\flag{\Gamma}_\infty$. These computations will be used later to relate the finiteness properties of $\BB_\Gamma$ to the topology of $\flag{\Gamma}$ and $\flag{\Gamma}_\infty$.

Recall that $\Delta$ is the graph associated to $\Gamma$ in \cref{section:brackets of cliques}, and that $\Delta_{\fin}$ is its full subgraph spanned by the finite vertices. By \cref{cor:contractible:links}, the links associated to cliques not contained in $\Delta_{\fin}$ are contractible. Thus only cliques in $\Delta_{\fin}$ require further consideration. For such cliques, the descending link has the same homology as the descending up-link, so it remains to determine the reduced homology of the ascending link.
\begin{lemma}\label[lemma]{lemma:homology:links}
    Let $C$ be a non-empty clique in $\Delta$ formed only by finite vertices.
    \begin{enumerate}
        \item If $|C|\ge 3$, then $\aslk[C]$ is simply connected. If $|C|=2$, then $\aslk[C]$ is simply connected if and only if $\lk(C,\flag{\Gamma})_\infty$ is non-empty. If $|C|=1$, then $\aslk[C]$ is simply connected if and only if $\lk(C,\flag{\Gamma})_\infty$ is non-empty and connected.
        \item For every integer $n$, the complex $\aslk[C]$ is homologically $(|C|+n)$-connected if and only if $\lk(C,\flag{\Gamma})_\infty$ is homologically $n$-connected. If $n\le -2$, both conditions hold automatically: homological $n$-connectivity is then a vacuous condition, and $\aslk[C]$ is always homologically $(|C|-2)$-connected. Moreover, for every $k\ge 0$,
        \[
            \widetilde{H}_k(\aslk[C])\cong
            \bigoplus_{I(C)}\widetilde{H}_{k-|C|}(\lk(C,\flag{\Gamma})_\infty),
        \]
        where reduced homology is understood in the augmented sense: $\widetilde{H}_{-1}(K)$ detects whether $K$ is empty, and $\widetilde{H}_j(K)=0$ for $j<-1$.
    \end{enumerate}
\end{lemma}
\begin{proof}
By \cref{cor:contractible:links}(3), for every clique $C$ there is a non-empty finite set of indices $I(C)$ such that
\[
    \aslk[C]\simeq\asulk[C]*\mvee_{I(C)}S^{|C|-1}.
\]
If $\asulk[C]$ is non-empty, this gives the chain of homotopy equivalences
\[
    \aslk[C]
    \simeq \Sigma\left(\asulk[C]\mwedge \mvee_{I(C)}S^{|C|-1}\right)
    \simeq \mvee_{I(C)}\Sigma^{|C|}\asulk[C].
\]
Hence $\aslk[C]$ is simply connected whenever $|C|\ge 2$ and $\asulk[C]$ is non-empty. If $|C|=1$ and $\asulk[C]$ is non-empty, then $\aslk[C]\simeq\mvee_{I(C)}\Sigma\asulk[C]$ with $I(C)\neq\emptyset$, which is simply connected if and only if $\asulk[C]$ is connected. If $\asulk[C]$ is empty, then
\[
    \aslk[C]\simeq\mvee_{I(C)}S^{|C|-1},
\]
which is simply connected exactly when $|C|\ge 3$.
By \cref{thm:description:as:des:links} and \cref{obs:delta:gamma:finite:links}, the complex $\asulk[C]$ is homotopy equivalent to the barycentric subdivision of $\lk(C,\flag{\Gamma})_\infty$ (see the end of the proof). In particular, $\asulk[C]$ is non-empty, respectively connected, if and only if $\lk(C,\flag{\Gamma})_\infty$ is, which proves the first assertion.

If $\asulk[C]$ is non-empty, the suspension equivalence above gives the following formula. If it is empty, the same formula follows from the wedge of spheres and the augmented homology convention. Thus, for every $k\ge 0$,
\[
    \widetilde{H}_k(\aslk[C])
    \cong\bigoplus_{I(C)}\widetilde{H}_{k-|C|}(\asulk[C]).
\]
By \cref{thm:description:as:des:links}, $\asulk[C]$ is homotopy equivalent to $\lk(C,\flag{\Delta})^b\cap (\flag{\Delta}_{>0})^b$. This complex is the barycentric subdivision of $\lk(C,\flag{\Delta})\cap\flag{\Delta}_{>0}$, which is isomorphic to $\lk(C,\flag{\Gamma})_\infty$ by \cref{obs:delta:gamma:finite:links}. Hence these complexes have the same reduced homology and homological connectivity.

Finally, $\aslk[C]$ is non-empty, and since $I(C)\neq\emptyset$, the formula shows that $\widetilde{H}_k(\aslk[C])=0$ for all $k\le |C|+n$ if and only if $\widetilde{H}_j(\lk(C,\flag{\Gamma})_\infty)=0$ for all $j\le n$; here indices $j\le -2$ contribute nothing by the augmented convention. This gives the equivalence in the second assertion for every integer $n$. If $n\le -2$, then $|C|+n\le|C|-2$ and every relevant group $\widetilde{H}_{k-|C|}(\lk(C,\flag{\Gamma})_\infty)$ has index $k-|C|\le -2$, so it vanishes; thus $\aslk[C]$ is homologically $(|C|-2)$-connected, and in particular homologically $(|C|+n)$-connected.
\end{proof}

\section{Proof of \cref{thm:newer:A}}\label{section:proof-main-a}

We first recall the Bestvina--Brady criterion relating the connectivity of Morse links to finiteness properties of kernels. Let $G$ act freely, properly discontinuously, cocompactly, and cellularly on a contractible affine cell complex $Y$, and suppose that the action preserves the affine structure.

\begin{theorem}\cite[Theorem 4.1]{BB}\label{BB link hom implies fin}
Let $A\leqslant\mathbb{R}$ be a nonzero infinite cyclic subgroup, let $\phi\colon G\to A$ be an epimorphism, and let $h:Y\to\mathbb{R}$ be a $\phi$-equivariant Morse function.
\begin{enumerate}
    \item Suppose that each ascending link and each descending link is homologically $n$-connected. Then $\ker\phi$ is of type $\fh_{n+1}(R)$.
    \item If all ascending and descending links are $R$-acyclic, then $\ker\phi$ is of type $\fh(R)$.
    \item If all ascending and descending links are simply connected, then $\ker\phi$ is finitely presented.
\end{enumerate}
\end{theorem}

In our setting, $G(\Gamma)$ acts on the Ruane--Witzel complex as above. However, the action is not free when $\Gamma$ has non-trivial finite vertices. The following consequence replaces freeness by finite stabilizers for virtually torsion-free groups.

\begin{corollary}\label{cor:BB:finitestabilizers}
    Let $G$ be a virtually torsion-free group acting with finite stabilizers, properly discontinuously, cocompactly, and cellularly on a contractible affine cell complex $Y$, preserving the affine structure. Let $A\leqslant\mathbb{R}$ be a nonzero infinite cyclic subgroup, let $\phi:G\to A$ be an epimorphism, and let $h:Y\to\mathbb{R}$ be a $\phi$-equivariant Morse function.
    \begin{enumerate}
    \item Suppose that each ascending link and each descending link is homologically $n$-connected. Then $\ker\phi$ is of type $\fp_{n+1}(R)$.
    \item If all ascending and descending links are simply connected, then $\ker\phi$ is finitely presented.
\end{enumerate}
\end{corollary}
\begin{proof}
Let $H\leqslant G$ be a torsion-free subgroup of finite index. Since the
$G$-action has finite stabilizers, its restriction to $H$ is free. Moreover,
the restricted action remains properly discontinuous, cocompact, and
cellular. Since $[G:H]$ is finite, $\phi(H)$ has finite index in $A$. Let
\[
\phi_H=\left.\phi\right|_H\colon H\longrightarrow \phi(H).
\]
The function $h$ is $\phi_H$-equivariant and has the same ascending and descending links as for the original $G$-action. Furthermore,
\[
\ker\phi_H=H\cap\ker\phi
\qquad\text{and}\qquad
[\ker\phi:\ker\phi_H]\leq [G:H]<\infty.
\]

For~(1), assume that every ascending and descending link of $h$ is
homologically $n$-connected. By \cref{BB link hom implies fin}(1),
$\ker\phi_H$ is of type $\fh_{n+1}(R)$, and hence of type
$\fp_{n+1}(R)$. Since type $\fp_{n+1}(R)$ is invariant under passage to
finite-index overgroups, $\ker\phi$ is also of type $\fp_{n+1}(R)$.

For~(2), \cref{BB link hom implies fin}(3) implies that $\ker\phi_H$ is
finitely presented. Finite presentability is invariant under passage to
finite-index overgroups, so $\ker\phi$ is finitely presented as well.
\end{proof}

To apply the corollary to $G(\Gamma)$, we need the following standard fact.

\begin{lemma}\label[lemma]{lemma:graph:product:vtf}
The graph product $G(\Gamma)$ is virtually torsion-free.
\end{lemma}
\begin{proof}
The vertex groups are finite or infinite cyclic, hence residually finite, so
$G(\Gamma)$ is residually finite by Green's theorem \cite{Green90}. Every
torsion element of $G(\Gamma)$ is conjugate into one of the finitely many
finite subgroups $[C_{\fin}]$. Since $G(\Gamma)$ is residually finite, there
is a finite quotient $\pi\colon G(\Gamma)\to Q$ that is injective on each of
them. If $g\in\ker\pi$ has finite order, then some conjugate of $g$ lies in
$\ker\pi\cap[C_{\fin}]=1$ for some $C$, because $\ker\pi$ is normal; hence
$g=1$. Thus $\ker\pi$ is a torsion-free subgroup of finite index.
\end{proof}

Set $\kappa=\dim\mathcal{S}+1$. We apply the corollary, which is possible by \cref{lemma:graph:product:vtf} and \cref{lemma:basic-properties-rw}(4), to the epimorphism $\kappa f\colon G(\Gamma)\to \kappa\mathbb{Z}$ and the Morse function $F\colon X\to\dbR$ from \cref{thm:F:fequivariant,thm:F:morsesimplicial}. Notice that $\ker(\kappa f)=\BB_\Gamma$. Recall that $\Delta=\Delta(\Gamma)$, and a \emph{clique of finite vertices} means a clique whose vertices all lie in $\Delta_{\fin}$.

\newpage

\begin{theorem}\label{thm:principal:connectivity implies finiteness}
Let $R$ be a nonzero commutative ring. The following statements hold.
\begin{enumerate}
    \item If $\flag{\Gamma}$ and $\flag{\Gamma}_\infty$ are homologically $n$-connected and, for every non-empty clique $C$ of finite vertices, $\lk(C,\flag{\Gamma})$ is homologically $n$-connected and $\lk(C,\flag{\Gamma})_\infty$ is homologically $(n-|C|)$-connected, then $\BB_\Gamma$ is of type $\fp_{n+1}(R)$.
    \item Suppose that $\flag{\Gamma}$ and $\flag{\Gamma}_\infty$ are simply connected and that $\lk(C,\flag{\Gamma})$ is simply connected for every non-empty clique $C$ of finite vertices. If $\lk(w,\flag{\Gamma})_\infty$ is non-empty and connected for every finite vertex $w$ and $\lk(C,\flag{\Gamma})_\infty$ is non-empty for every clique $C$ of finite vertices with $|C|=2$, then $\BB_\Gamma$ is finitely presented.
\end{enumerate}
\end{theorem}
\begin{proof}
By \cref{cor:contractible:links}, only the empty clique and the non-empty cliques of finite vertices need to be considered. Their descending links are controlled by $\flag{\Gamma}$ and $\lk(C,\flag{\Gamma})$, and their ascending links by $\flag{\Gamma}_\infty$ and \cref{lemma:homology:links}. For~(1), \cref{lemma:homology:links}(2) applied with $n-|C|$ in place of $n$ shows that $\aslk[C]$ is homologically $n$-connected; this includes the case $n-|C|\le -2$, where the hypothesis on $\lk(C,\flag{\Gamma})_\infty$ is vacuous. For finite presentation, \cref{lemma:homology:links}(1) shows that $\aslk[C]$ is simply connected: for $|C|=1$ by the connectivity of $\lk(w,\flag{\Gamma})_\infty$, for $|C|=2$ by the non-emptiness assumption, and for $|C|\ge3$ with no further assumption. The result follows from the preceding corollary.
\end{proof}

\section{Homology of level sets}\label{section:homology-level-sets}

Given $j\in\dbZ$, define $X_j$, $X_{\le j}$, respectively $X_{\ge j}$, as
$F^{-1}(j)$, $F^{-1}(-\infty,j]$, respectively $F^{-1}[j,\infty)$.
We write $CY$ for the cone on a space $Y$.

Throughout this section, $R$ denotes a nonzero commutative ring. We use
augmented reduced homology, with
\[
    \widetilde H_{-1}(\emptyset;R)=R,
    \qquad
    \widetilde H_k(L;R)=0\quad\text{for }k<-1.
\]
For $q\in\dbZ$, we say that a simplicial complex $L$ is
\emph{strictly $(R,q)$-acyclic} if
\[
    \widetilde H_k(L;R)=0
    \qquad\text{for }k=q+2\text{ and for }-1\le k\le q,
\]
and
\[
    \widetilde H_{q+1}(L;R)\neq0.
\]
The second range is understood to be empty when $q<-1$. We say that $L$
satisfies the \emph{$(R,q)$-vanishing condition} if it satisfies only the
vanishing condition above.

The next proposition is an analogue of \cite[Proposition~7.1]{BB}. 

\begin{proposition}\label[proposition]{prop:homology:level:sets:complex}
Let $R$ be a nonzero commutative ring, let $n\ge0$, and let $j\in\dbZ$.
The following statements hold.
\begin{enumerate}
    \item Suppose that $\deslk[C]$ satisfies the $(R,n-1)$-vanishing
    condition for every clique $C$ which is either empty or consists of finite
    vertices. Then
    \[
        H_k(X,X_{\le j};R)=0,
        \qquad 0\le k\le n.
    \]
    If, in addition, $\deslk[C_0]$ is strictly $(R,n-1)$-acyclic for at least
    one such clique $C_0$, then $H_{n+1}(X,X_{\le j};R)$ is not finitely
    generated as an $R[\BB_\Gamma]$-module.

    \item Suppose that $\aslk[C]$ satisfies the $(R,n-1)$-vanishing
    condition for every clique $C$ which is either empty or consists of finite
    vertices. Then
    \[
        H_k(X,X_{\ge j};R)=0,
        \qquad 0\le k\le n.
    \]
    If, in addition, $\aslk[C_0]$ is strictly $(R,n-1)$-acyclic for at least
    one such clique $C_0$, then $H_{n+1}(X,X_{\ge j};R)$ is not finitely
    generated as an $R[\BB_\Gamma]$-module.
\end{enumerate}
\end{proposition}

\begin{proof}
We first prove (1). For $m\in\dbN$, set $Y_m=X_{\le j+m}$, so that
$Y_0=X_{\le j}$ and $X=\bigcup_{m\ge0}Y_m$. The Morse lemma gives
\begin{equation}\label{eq:relative:morse:level:complex}
    H_k(Y_m,Y_{m-1};R)
    \cong
    \bigoplus_{F(v)=j+m}H_k(C\deslk(v),\deslk(v);R)
    \cong
    \bigoplus_{F(v)=j+m}\widetilde H_{k-1}(\deslk(v);R).
\end{equation}
If the type of a vertex contains an infinite vertex, its ascending and
descending links are contractible. Thus only vertices whose type is the
empty clique or a clique of finite vertices can contribute to these relative
homology groups.

The vanishing hypothesis and \eqref{eq:relative:morse:level:complex} give
\[
    H_k(Y_m,Y_{m-1};R)=0,
    \qquad 0\le k\le n,
\]
and
\[
    H_{n+2}(Y_m,Y_{m-1};R)=0.
\]
Using the long exact sequence of the triple $(Y_m,Y_{m-1},Y_0)$, induction
on $m$ gives
\begin{equation}\label{eq:relative:vanishing:level:complex}
    H_k(Y_m,Y_0;R)=0,
    \qquad 0\le k\le n.
\end{equation}
Passing to the filtered colimit gives
\[
    H_k(X,X_{\le j};R)=0,
    \qquad 0\le k\le n.
\]

The same long exact sequence, together with the vanishing in degrees $n$ and
$n+2$, gives a short exact sequence of $R[\BB_\Gamma]$-modules
\begin{equation}\label{eq:short:exact:relative:levels:complex}
    0\longrightarrow H_{n+1}(Y_{m-1},Y_0;R)
    \longrightarrow H_{n+1}(Y_m,Y_0;R)
    \longrightarrow H_{n+1}(Y_m,Y_{m-1};R)
    \longrightarrow0.
\end{equation}
By \eqref{eq:relative:morse:level:complex}, the last term is
\begin{equation}\label{eq:relative:piece:complex}
    \bigoplus_{F(v)=j+m}\widetilde H_n(\deslk(v);R).
\end{equation}

Assume now that $\deslk[C_0]$ is strictly $(R,n-1)$-acyclic. Let $u$ be an
infinite standard generator of $G(\Gamma)$ and set $\kappa=\dim\mathcal S+1$.
Equivariance gives
\[
    F(u^\ell[C_0])=F([C_0])+\kappa\ell.
\]
Hence vertices of type $[C_0]$ occur at infinitely many levels whose heights
tend to $+\infty$. Since
$\widetilde H_n(\deslk[C_0];R)\neq0$, the module in
\eqref{eq:relative:piece:complex} is nonzero for arbitrarily large $m$.

The first map in \eqref{eq:short:exact:relative:levels:complex} is injective,
and its image is proper whenever the last term is nonzero. After identifying
the modules by these injections, they form an increasing sequence with
strict inclusions for arbitrarily large $m$. Since homology commutes with
filtered colimits,
\[
    H_{n+1}(X,X_{\le j};R)
    \cong
    \varinjlim_m H_{n+1}(Y_m,Y_0;R).
\]
If this module were generated by finitely many elements, all of them would
belong to the image of a single stage. That stage would then contain the
whole colimit, contradicting a later strict inclusion. Therefore
$H_{n+1}(X,X_{\le j};R)$ is not finitely generated over
$R[\BB_\Gamma]$.

The proof of (2) is the same, using the exhaustion $Y_m=X_{\ge j-m}$ and
coning off ascending links. The nonzero relative terms occur at infinitely
many levels whose heights tend to $-\infty$.
\end{proof}

The preceding proposition gives the following analogue of
\cite[Corollary~7.2]{BB}.

\begin{corollary}\label[corollary]{cor:homology:level:set:infinitely:generated}
Let $R$ be a nonzero commutative ring, let $n\ge0$, and let $j\in\dbZ$.
Suppose that $\deslk[C]$ and $\aslk[C]$
satisfy the $(R,n-1)$-vanishing condition for every clique $C$ which is
either empty or consists of finite vertices. Suppose moreover that, for at
least one such clique $C_0$, either $\deslk[C_0]$ or $\aslk[C_0]$ is strictly
$(R,n-1)$-acyclic. Then
\[
    \widetilde H_k(X_j;R)=0,
    \qquad 0\le k\le n-1,
\]
and $\widetilde H_n(X_j;R)$ is not finitely generated as an
$R[\BB_\Gamma]$-module.
\end{corollary}

\begin{proof}
By \cref{prop:homology:level:sets:complex},
\[
    H_k(X,X_{\le j};R)=0
    \quad\text{and}\quad
    H_k(X,X_{\ge j};R)=0,
    \qquad 0\le k\le n.
\]
Since $X$ is contractible, the long exact sequences of the two pairs give
\[
    \widetilde H_k(X_{\le j};R)=0
    \quad\text{and}\quad
    \widetilde H_k(X_{\ge j};R)=0,
    \qquad 0\le k\le n-1.
\]
The reduced Mayer--Vietoris sequence for
$X=X_{\le j}\cup X_{\ge j}$ gives
\[
    \widetilde H_k(X_j;R)=0,
    \qquad 0\le k\le n-1.
\]
In degree $n$, the long exact sequences of the pairs and Mayer--Vietoris give
\[
    \widetilde H_n(X_j;R)
    \cong
    \widetilde H_n(X_{\le j};R)
    \oplus
    \widetilde H_n(X_{\ge j};R).
\]
At least one summand is not finitely generated by
\cref{prop:homology:level:sets:complex}. This argument also applies when
$n=0$, using augmented reduced homology.
\end{proof}

The preceding corollary can be combined with the explicit descriptions of
ascending and descending links from \cref{thm:description:as:des:links} and
\cref{lemma:homology:links}. This gives the following partial converse of
\cref{thm:principal:connectivity implies finiteness}.

\begin{theorem}\label{thm:partial:converse:principal:links}
Let $R$ be a nonzero commutative ring and let $n\ge0$. Suppose that:
\begin{enumerate}
    \item $\flag{\Gamma}$ and, for every non-empty clique $C$ of finite
    vertices, $\lk(C,\flag{\Gamma})$ satisfy the
    $(R,n-1)$-vanishing condition.

    \item $\flag{\Gamma}_\infty$ satisfies the $(R,n-1)$-vanishing condition
    and, for every non-empty clique $C$ of finite vertices,
    $\lk(C,\flag{\Gamma})_\infty$ satisfies the
    $(R,n-|C|-1)$-vanishing condition.

    \item At least one of the following complexes is strictly acyclic in the
    indicated degree:
    \begin{enumerate}
        \item $\flag{\Gamma}$ is strictly $(R,n-1)$-acyclic;

        \item $\lk(C_0,\flag{\Gamma})$ is strictly $(R,n-1)$-acyclic for
        some non-empty clique $C_0$ of finite vertices;

        \item $\flag{\Gamma}_\infty$ is strictly $(R,n-1)$-acyclic;

        \item $\lk(C_0,\flag{\Gamma})_\infty$ is strictly
        $(R,n-|C_0|-1)$-acyclic for some non-empty clique $C_0$ of finite
        vertices.
    \end{enumerate}
\end{enumerate}
Then $\BB_\Gamma$ is not of type $\fp_{n+1}(R)$.
\end{theorem}

\begin{proof}
The case $C=\emptyset$ for descending links is the condition on
$\flag{\Gamma}$: the complex $\deslk[\emptyset]$ is naturally identified
with the barycentric subdivision of $\flag{\Gamma}$. If
$C\subseteq\Delta_{\fin}$ is non-empty, then
\[
    \deslk[C]\cong\desulk[C]
    \cong\lk(C,\flag{\Gamma})^b.
\]

The case $C=\emptyset$ for ascending links identifies
$\aslk[\emptyset]$ with the barycentric subdivision of
$\flag{\Gamma}_\infty$. For a non-empty clique
$C\subseteq\Delta_{\fin}$, \cref{lemma:homology:links} gives
\[
    \widetilde H_k(\aslk[C];R)
    \cong
    \bigoplus_{I(C)}
    \widetilde H_{k-|C|}(\lk(C,\flag{\Gamma})_\infty;R)
\]
for every $k\ge0$. This formula remains valid when
$\lk(C,\flag{\Gamma})_\infty$ is empty by the augmented homology convention.
Since $I(C)$ is non-empty, the hypotheses imply that all ascending and
descending links satisfy the $(R,n-1)$-vanishing condition and that at least
one of them is strictly $(R,n-1)$-acyclic. Hence
\cref{cor:homology:level:set:infinitely:generated} gives
\[
    \widetilde H_i(X_j;R)=0,
    \qquad 0\le i\le n-1,
\]
and $\widetilde H_n(X_j;R)$ is not finitely generated over
$R[\BB_\Gamma]$.

Let $H=\BB_\Gamma$. By \cref{lemma:graph:product:vtf}, $G(\Gamma)$ is
virtually torsion-free, and hence so is $H$; choose a torsion-free finite-index subgroup
$H_0\le H$. The $H$-action on $X_j$ is cocompact by the cocompactness of the
$G(\Gamma)$-action on $X$ and the $\kappa f$-equivariance of $F$. The intersections
of the simplices of $X$ with $X_j$ therefore give an $H_0$-free, cocompact
affine cell structure on $X_j$. If
$\widetilde H_n(X_j;R)$ were finitely generated over $R[H_0]$, then it would
also be finitely generated over $R[H]$. Thus it is not finitely generated
over $R[H_0]$.

Applying \cite[Lemma~3.8]{BB} to the action of $H_0$ on $X_j$ shows that
$H_0$ is not of type $\fp_{n+1}(R)$. Since type $\fp_{n+1}(R)$ is invariant
under passage to finite-index subgroups and overgroups, the same holds for
$H$.
\end{proof}

The first two conditions of
\cref{thm:partial:converse:principal:links} imply, by
\cref{thm:principal:connectivity implies finiteness}, that $\BB_\Gamma$ is of
type $\fp_n(R)$. Thus the theorem produces groups of type $\fp_n(R)$ but not
of type $\fp_{n+1}(R)$ over arbitrary nonzero commutative coefficient rings.

\section{Bredon finiteness properties}\label{section:bredon-finiteness}

Throughout this section, we work with finite labeled graphs whose vertex
groups are infinite cyclic or finite abelian. This setting is convenient for
describing normalizers of finite subgroups as graph products.

Fix a nonzero commutative ring $R$ and let $n\geq 0$. A standard Bredon
finiteness criterion (see \cite{KMPN09}) says that
$\BB_\Gamma$ is of type $\underline{\fp}_n(R)$ if and only if it has
finitely many conjugacy classes of finite subgroups and every finite
subgroup $Q\leqslant\BB_\Gamma$ has a normalizer of type $\fp_n(R)$.

We study these two conditions in turn. We first characterize when
$\BB_\Gamma$ has finitely many conjugacy classes of finite subgroups, and
then use the normalizer computation to study when the normalizers have type
$\fp_n(R)$.

\subsection{Conjugacy classes of finite subgroups}\label{section:bredon-conjugacy-classes}
\begin{theorem}\label[theorem]{thm:finite-conjugacy-classes-BB}
Suppose that $\Gamma$ is finite and that its vertex groups are infinite
cyclic or finite abelian. Then $\BB_\Gamma$ has finitely many conjugacy
classes of finite subgroups if and only if
\(
    \lk(C,\flag{\Gamma})_\infty\neq\varnothing
\)
for every non-empty clique $C$ of finite vertices. If this condition holds,
each $G(\Gamma)$-conjugacy class of finite subgroups gives exactly one
$\BB_\Gamma$-conjugacy class.
\end{theorem}

\begin{proof}
Every finite subgroup of $G(\Gamma)$ fixes a point in its Ruane--Witzel complex
$X$. It is therefore conjugate into the stabilizer of a simplex, and hence
into $[C]$ for some clique $C$ of finite vertices. Since $\Gamma$ is finite,
there are only finitely many such cliques, and each $[C]$ is finite. Thus
$G(\Gamma)$ has only finitely many conjugacy classes of finite subgroups.
Every finite subgroup lies in $\BB_\Gamma$, since its image under the
homomorphism $f:G(\Gamma)\to\dbZ$ is a finite subgroup of $\dbZ$ and is
therefore trivial.

Fix a finite subgroup $Q\leqslant\BB_\Gamma$. The $\BB_\Gamma$-orbits in
its $G(\Gamma)$-conjugacy class are parametrized by
\[
    \BB_\Gamma\backslash G(\Gamma)/N_{G(\Gamma)}Q
    \cong \dbZ/f\bigl(N_{G(\Gamma)}Q\bigr).
\]
The double coset of $g$ corresponds to the $\BB_\Gamma$-conjugacy class of $gQg^{-1}$; the displayed bijection sends it to $f(g)+f(N_{G(\Gamma)}Q)$.
For $Q=1$ this quotient is trivial. Otherwise, conjugate $Q$ by an element of $G(\Gamma)$ so that its
support is a non-empty clique $C$ of finite vertices. By
\cref{prop:normalizer:finite:revised},
\[
    N_{G(\Gamma)}Q=G(\St_\Gamma(C)).
\]
As $C$ consists of finite vertices, the image of this normalizer under
$f$ is $\dbZ$ precisely when
$\lk(C,\flag{\Gamma})_\infty$ is non-empty; otherwise its image is zero.
Thus the $G(\Gamma)$-conjugacy class splits into one or infinitely many
$\BB_\Gamma$-conjugacy classes, respectively. In particular, if the link
condition fails for $C$, the finite subgroup $[C]$ has infinitely many
$\BB_\Gamma$-conjugacy classes among its $G(\Gamma)$-conjugates. Together
with the finiteness of the number of $G(\Gamma)$-conjugacy classes, this
proves the assertion.
\end{proof}

Consequently, under the hypotheses of
\cref{thm:finite-conjugacy-classes-BB}, the finite-conjugacy-class
condition in the Bredon criterion holds precisely when every non-empty
clique of finite vertices has an infinite common neighbor.

\subsection{Finiteness properties of normalizers}\label{section:bredon-normalizers}

We now study the second condition in the Bredon criterion. We first compute
the normalizers algebraically and then apply the finiteness criteria for
Bestvina--Brady kernels.

For $g\in G(\Gamma)$, choose a minimal length factorization
$g=g_1\cdots g_m$ with each $g_i$ in a vertex group. The \emph{support} of
$g$ is the set of vertices $v$ for which $g_i\in G_v\setminus\{1\}$ for
some $i$. For a subgroup $Q\leqslant G(\Gamma)$, let
$\supp_\Gamma(Q)$ be the full subgraph generated by the supports of its
elements. We recall the following results from \cite{AM11}.
\begin{proposition}\label[proposition]{prop:normalizer:finite:revised}
Let $\Gamma$ be a labeled graph whose vertex groups are infinite cyclic or
finite abelian. Let $1\neq Q\leqslant G(\Gamma)$ be finite and suppose that
$\supp_\Gamma(Q)$ spans a clique $\calC$ of finite vertices. Then
\[
    N_{G(\Gamma)}Q=G(\St_\Gamma(\calC)).
\]
\end{proposition}

\begin{proof}
The parabolic closure of $Q$, i.e. the minimal subgroup of $G(\Gamma)$ of the form $G(\Gamma')$ with $\Gamma'$ a full subgraph of $\Gamma$, is the standard parabolic subgroup $[\calC]$.
By \cite[Lemma~3.12]{AM11}, every element normalizing $Q$ normalizes its
parabolic closure. Hence
\[
    N_{G(\Gamma)}Q
    \leqslant N_{G(\Gamma)}[\calC]
    =G(\St_\Gamma(\calC)),
\]
where the equality follows from \cite[Proposition~3.13]{AM11}.

Conversely, every vertex group corresponding to a vertex of
$\lk(\calC,\flag{\Gamma})$ centralizes $Q$. Every vertex group corresponding
to a vertex of $\calC$ normalizes $Q$, because
$Q\leqslant[\calC]$ and $[\calC]$ is abelian. These vertex groups generate
$G(\St_\Gamma(\calC))$, proving the reverse inclusion.
\end{proof}

\begin{theorem}\label{thm:normalizer:BB:revised}
Let $\Gamma$ be as above and let $Q$ be a finite subgroup of $\BB_\Gamma$.
If $Q=1$, then $N_{\BB_\Gamma}Q=\BB_\Gamma$. If $Q\neq1$, then there exist
$g\in G(\Gamma)$ and a non-empty clique $\calC$ of finite vertices such that
\[
    N_{\BB_\Gamma}Q
    =g^{-1}\BB_{\St_\Gamma(\calC)}g.
\]
In particular,
\[
    N_{\BB_\Gamma}Q\cong\BB_{\St_\Gamma(\calC)}.
\]
\end{theorem}

\begin{proof}
Suppose that $Q\neq1$. A finite subgroup of $G(\Gamma)$ is conjugate into a
standard parabolic subgroup $[C]$ generated by a clique of finite vertices; indeed this follows from the fact that every finite subgroup of $G(\Gamma)$ fixes a vertex in its Ruane--Witzel complex and the stabilizer of such a vertex is conjugate to a standard finite parabolic subgroup.
Choose $g\in G(\Gamma)$ such that $gQg^{-1}\leqslant[C]$, and let
$\calC=\supp_\Gamma(gQg^{-1})$. By
\cref{prop:normalizer:finite:revised},
\[
    N_{G(\Gamma)}(gQg^{-1})=G(\St_\Gamma(\calC)).
\]
Since $\BB_\Gamma$ is normal in $G(\Gamma)$,
\begin{align*}
    N_{\BB_\Gamma}Q
    &=N_{G(\Gamma)}Q\cap\BB_\Gamma\\
    &=g^{-1}\bigl(G(\St_\Gamma(\calC))\cap\BB_\Gamma\bigr)g\\
    &=g^{-1}\BB_{\St_\Gamma(\calC)}g.
\end{align*}
\end{proof}

The preceding theorem identifies each non-trivial finite-subgroup normalizer with a Bestvina--Brady kernel on a star. For the next two statements, choose
$g$ and $\calC$ as in \cref{thm:normalizer:BB:revised} and put
$\Lambda=\St_\Gamma(\calC)$.

\begin{corollary}\label[corollary]{cor:normalizer:positive:revised}
Let $1\neq Q\leqslant\BB_\Gamma$ be finite. The following statements hold.
\begin{enumerate}
    \item If $\flag{\Lambda}$ and $\flag{\Lambda}_\infty$ are homologically
    $n$-connected and, for every non-empty clique $C$ of finite vertices of
    $\Lambda$, the complexes $\lk(C,\flag{\Lambda})$ and
    $\lk(C,\flag{\Lambda})_\infty$ are respectively homologically
    $n$-connected and homologically $(n-|C|)$-connected, then
    $N_{\BB_\Gamma}Q$ is of type $\fp_{n+1}(R)$.

    \item If $\flag{\Lambda}$ and $\flag{\Lambda}_\infty$ are simply
    connected, and $\lk(C,\flag{\Lambda})$ is simply connected for every
    non-empty clique $C$ of finite vertices of $\Lambda$, while
    $\lk(w,\flag{\Lambda})_\infty$ is non-empty and connected for every finite vertex
    $w$ of $\Lambda$ and $\lk(C,\flag{\Lambda})_\infty$ is non-empty for
    every clique $C$ of finite vertices of $\Lambda$ with $|C|=2$, then
    $N_{\BB_\Gamma}Q$ is finitely presented.
\end{enumerate}
\end{corollary}

\begin{proof}
By \cref{thm:normalizer:BB:revised}, the normalizer is isomorphic to
$\BB_\Lambda$. Both assertions follow directly from
\cref{thm:principal:connectivity implies finiteness} applied to $\Lambda$.
\end{proof}

\begin{corollary}\label[corollary]{cor:normalizer:negative:revised}
Let $R$ be a nonzero commutative ring and let $n\geq0$.
Suppose that:
\begin{enumerate}
    \item $\flag{\Lambda}$ and, for every non-empty clique $C$ of finite
    vertices of $\Lambda$, $\lk(C,\flag{\Lambda})$ satisfy the
    $(R,n-1)$-vanishing condition.

    \item $\flag{\Lambda}_\infty$ satisfies the $(R,n-1)$-vanishing
    condition and, for every non-empty clique $C$ of finite vertices of
    $\Lambda$, $\lk(C,\flag{\Lambda})_\infty$ satisfies the
    $(R,n-|C|-1)$-vanishing condition.

    \item At least one of the following complexes is strictly acyclic in
    the indicated degree:
    \begin{enumerate}
        \item $\flag{\Lambda}$ is strictly $(R,n-1)$-acyclic;

        \item $\lk(C_0,\flag{\Lambda})$ is strictly $(R,n-1)$-acyclic for
        some non-empty clique $C_0$ of finite vertices of $\Lambda$;

        \item $\flag{\Lambda}_\infty$ is strictly $(R,n-1)$-acyclic;

        \item $\lk(C_0,\flag{\Lambda})_\infty$ is strictly
        $(R,n-|C_0|-1)$-acyclic for some non-empty clique $C_0$ of finite
        vertices of $\Lambda$.
    \end{enumerate}
\end{enumerate}
Then $N_{\BB_\Gamma}Q$ is not of type $\fp_{n+1}(R)$.
\end{corollary}

\begin{proof}
The group $N_{\BB_\Gamma}Q$ is isomorphic to $\BB_\Lambda$, and the result
is exactly \cref{thm:partial:converse:principal:links} applied to $\Lambda$.
\end{proof}

We now specialize to graphs in which no two finite vertices are adjacent.
This gives a convenient source of normalizers with precisely prescribed
homological finiteness properties. We say that a group is \emph{strictly of type $\fp_n(R)$} if it is of type
$\fp_n(R)$ but not of type $\fp_{n+1}(R)$.

\begin{corollary}\label[corollary]{cor:normalizer:isolated:revised} 
Let $v$ be a finite vertex of $\Gamma$, assume that no two finite vertices of
$\Gamma$ are adjacent, and let $1\neq Q\leqslant G_v$. Set
\[
    L=\lk(v,\flag{\Gamma}),
\]
and let $R$ be a nonzero commutative ring.
Then:
\begin{enumerate}
    \item if $L$ is homologically $n$-connected, then
    $N_{\BB_\Gamma}Q$ is of type $\fp_{n+1}(R)$;
    \item if $L$ is simply connected, then $N_{\BB_\Gamma}Q$ is finitely
    presented.
    \item if $L$ is strictly $(R,n)$-acyclic, then $N_{\BB_\Gamma}Q$ is strictly of type $\fp_{n+1}(R)$.
\end{enumerate}
\end{corollary}

\begin{proof}
Here $\supp_\Gamma(Q)=\{v\}$, so
\[
    N_{\BB_\Gamma}Q\cong\BB_{\St_\Gamma(v)}.
\]
Since no finite vertex is adjacent to $v$, every vertex of $L$ is infinite.
Moreover, $\flag{\St_\Gamma(v)}$ is the cone $v*L$, its infinite subcomplex
is $L$, and the only non-empty clique of finite vertices is $\{v\}$, whose
link and infinite link are both $L$. Thus the first two items follow from
\cref{cor:normalizer:positive:revised}.

Strict $(R,n)$-acyclicity implies that $L$ is homologically $n$-connected,
so the first item gives type $\fp_{n+1}(R)$.
For the last item, note that $\flag{\St_\Gamma(v)}=v*L$ and $v$ is
the only finite vertex of $\St_\Gamma(v)$, so
$\BB_{\St_\Gamma(v)}\cong G_v\times\BB_L$ with $G_v$ finite, where $\BB_L$ is
the Bestvina--Brady group of the graph product of infinite cyclic groups on
$L$. Since $G_v$ is finite, $G_v\times\BB_L$ is of type $\fp_{m}(R)$ if and
only if $\BB_L$ is, because $\BB_L$ is a finite-index subgroup of it. Strict
$(R,n)$-acyclicity of $L$ means in particular that
$\widetilde H_k(L;R)=0$ for $-1\le k\le n$ and
$\widetilde H_{n+1}(L;R)\neq0$, so the classical Bestvina--Brady theorem
applied to $\BB_L$ shows that it is of type $\fp_{n+1}(R)$ but not of type
$\fp_{n+2}(R)$. Hence the normalizer is strictly of type $\fp_{n+1}(R)$.
\end{proof}

Under the Bredon characterization stated above, the preceding results give
the following criterion.

\begin{theorem}\label{thm:bredon-fp-criterion}
Let $R$ be a nonzero commutative ring, let $n\geq 0$, and suppose that
$\Gamma$ is finite with vertex groups infinite cyclic or finite abelian.
Then $\BB_\Gamma$ is of type $\underline{\fp}_n(R)$ if and only if:
\begin{enumerate}
    \item $\lk(C,\flag{\Gamma})_\infty\neq\varnothing$ for every
    non-empty clique $C$ of finite vertices; and
    \item $\BB_\Gamma$ is of type $\fp_n(R)$, and
    $\BB_{\St_\Gamma(C)}$ is of type $\fp_n(R)$ for every non-empty clique
    $C$ of finite vertices.
\end{enumerate}
\end{theorem}

\begin{proof}
By the Bredon characterization stated above,
$\BB_\Gamma$ is of type $\underline{\fp}_n(R)$ precisely when it has
finitely many conjugacy classes of finite subgroups and every finite
subgroup $Q\leqslant\BB_\Gamma$ has a normalizer of type $\fp_n(R)$.
By \cref{thm:finite-conjugacy-classes-BB}, the first condition is equivalent
to item (1). By \cref{thm:normalizer:BB:revised}, the trivial subgroup has
normalizer $\BB_\Gamma$, while the normalizer of every non-trivial finite
subgroup is conjugate to $\BB_{\St_\Gamma(C)}$ for some non-empty clique $C$
of finite vertices. Conversely, each such clique occurs as the support of
the finite subgroup $[C]\leqslant\BB_\Gamma$. Thus the normalizer condition
is equivalent to item (2), proving the result.
\end{proof}

\section{Two examples}\label{section:example}

\subsection{First example}

We apply our results to construct Bestvina--Brady kernels in which one can choose the strict finiteness properties of the normalizers of the maximal finite subgroups. In particular, this gives a simple construction of groups which are strictly of type $\underline{\fp}_n(R)$ with a prescribed number of conjugacy classes of maximal finite subgroups.

\begin{proposition}\label[proposition]{prop:example}
Let $R$ be a nonzero commutative ring, and let $m\ge 1$, $n\ge 1$ and $n_1,\dots,n_m\ge n$ be integers. There is a labeled graph $\Gamma$ such that $\BB_\Gamma$ is of type $\fp_n(R)$, and of type $\textup{F}_n$ if $n\ge 2$, and has exactly $m$ conjugacy classes of maximal finite subgroups, with representatives $Q_1,\dots,Q_m$ such that $N_{\BB_\Gamma}Q_i$ is strictly of type $\fp_{n_i}(R)$ (and strictly of type $\textup{F}_{n_i}$ if $n_i\ge 2$). If moreover $n_i=n$ for all $i$, then $\BB_\Gamma$ is strictly of type $\underline{\fp}_n(R)$.
\end{proposition}

\begin{proof}
Let $Y$ be $S^{n+1}$ with the interiors of $m$ disjoint $(n+1)$-discs removed. For each $i$, glue an $(n_i+1)$-disc $E_i$ to $Y$ by identifying the $i$-th boundary sphere of $Y$ with a PL $n$-sphere in $\partial E_i\cong S^{n_i}$, and call the result $X$. Triangulate $Y\cup\bigcup_i\partial E_i$ with each $\partial E_i$ a subcomplex, pass to the barycentric subdivision (which is flag, with each $\partial E_i$ full), and cone off each $\partial E_i$ with a new vertex $v_i$; this gives a flag triangulation of $X$. Let $\Gamma$ be its $1$-skeleton, with each $v_i$ labeled by a non-trivial finite abelian group $G_{v_i}$ and all other vertices by $\mathbb{Z}$.

Then $\flag{\Gamma}=X\simeq S^{n+1}$, the $v_i$ are pairwise non-adjacent, $\lk(v_i,\flag{\Gamma})=\lk(v_i,\flag{\Gamma})_\infty=\partial E_i\cong S^{n_i}$, and $\flag{\Gamma}_\infty=Y\cup\bigcup_i\partial E_i$, which is $(n-1)$-connected since it is obtained from $Y\simeq\bigvee_{m-1}S^n$ by gluing the spheres $\partial E_i$ along $n$-spheres. Hence \cref{thm:principal:connectivity implies finiteness} (with $n-1$ in place of $n$) shows that $\BB_\Gamma$ is of type $\fp_n(R)$, and in particular of type $\fp_n(\dbZ)$. If $n\ge 2$, then $\flag{\Gamma}$, $\flag{\Gamma}_\infty$ and the links $\lk(v_i,\flag{\Gamma})=\lk(v_i,\flag{\Gamma})_\infty\cong S^{n_i}$ are simply connected, and there are no cliques of finite vertices with two elements, so $\BB_\Gamma$ is finitely presented by \cref{thm:principal:connectivity implies finiteness}(2). Since a finitely presented group of type $\fp_n(\dbZ)$ is of type $\textup{F}_n$ \cite[Chapter~VIII, \S7]{Br94}, the group $\BB_\Gamma$ is of type $\textup{F}_n$ if $n\ge 2$.

As the finite vertices are pairwise non-adjacent, every finite subgroup of $G(\Gamma)$ is conjugate in $G(\Gamma)$ to a subgroup of some $G_{v_i}$, and $G_{v_i}$ is not conjugate to $G_{v_j}$ for $i\neq j$, because the retraction $G(\Gamma)\to G_{v_i}$ kills $G_{v_j}$. Since the non-empty cliques of finite vertices are the singletons $\{v_i\}$ and $\lk(v_i,\flag{\Gamma})_\infty\neq\varnothing$, \cref{thm:finite-conjugacy-classes-BB} shows that each $G(\Gamma)$-conjugacy class of finite subgroups gives exactly one $\BB_\Gamma$-conjugacy class; so the same statements hold with conjugacy in $\BB_\Gamma$. Hence $Q_i:=G_{v_i}$, $i=1,\dots,m$, represent the conjugacy classes of maximal finite subgroups, and $N_{\BB_\Gamma}Q_i$ is strictly of type $\fp_{n_i}(R)$ (and strictly of type $\textup{F}_{n_i}$ if $n_i\ge 2$) by \cref{cor:normalizer:isolated:revised}, since $\lk(v_i,\flag{\Gamma})\cong S^{n_i}$. For the type $\textup{F}_{n_i}$, recall from the proof of \cref{cor:normalizer:isolated:revised} that $N_{\BB_\Gamma}Q_i\cong G_{v_i}\times\BB_L$ with $L=\lk(v_i,\flag{\Gamma})\cong S^{n_i}$, and $\BB_L$ has finite index in it. Since $L$ is $(n_i-1)$-connected, the homotopical form of the Bestvina--Brady theorem \cite{BB} shows that $\BB_L$, and hence $N_{\BB_\Gamma}Q_i$, is of type $\textup{F}_{n_i}$. It is not of type $\textup{F}_{n_i+1}$, because type $\textup{F}_{n_i+1}$ implies type $\fp_{n_i+1}(R)$. Finally, let $n_i=n$ for all $i$. Every $1\neq Q\leqslant Q_i$ has support $\{v_i\}$, so $N_{\BB_\Gamma}Q=N_{\BB_\Gamma}Q_i$ by \cref{thm:normalizer:BB:revised}. Thus \cref{thm:bredon-fp-criterion} applies, as $\lk(v_i,\flag{\Gamma})_\infty\neq\varnothing$, so $\BB_\Gamma$ is of type $\underline{\fp}_n(R)$, but not of type $\underline{\fp}_{n+1}(R)$ because the $N_{\BB_\Gamma}Q_i$ are not of type $\fp_{n+1}(R)$.
\end{proof}

\subsection{Second example}\label{section:new-examples}

In contrast with \cref{section:example}, we now let all finite vertices span a single clique.

Let $k\geq1$ and let $\Gamma_k$ be the labeled graph with vertices $w_1,\dots,w_k,a_1,\dots,a_k$, in which two distinct vertices are adjacent unless they form one of the pairs $\{w_i,a_i\}$. Put $G_{w_i}=\langle s_i\rangle\cong\mathbb{Z}/2\mathbb{Z}$ and $G_{a_i}=\langle t_i\rangle\cong\mathbb{Z}$. Thus $\Gamma_k$ is the join of the $k$ edgeless graphs $\{w_i,a_i\}$, and
\[
    G(\Gamma_k)=\prod_{i=1}^k\bigl(\langle s_i\rangle*\langle t_i\rangle\bigr)\cong\prod_{i=1}^k\bigl(\mathbb{Z}/2\mathbb{Z}*\mathbb{Z}\bigr),
    \qquad
    f(s_i)=0,\quad f(t_i)=1.
\]
The flag complex $\flag{\Gamma_k}=\ast_{i=1}^k\{w_i,a_i\}$ is the boundary of the $k$-dimensional cross-polytope, so $\flag{\Gamma_k}\cong S^{k-1}$, while $\flag{\Gamma_k}_\infty$ is the $(k-1)$-simplex spanned by $a_1,\dots,a_k$. Let $W=\{w_1,\dots,w_k\}$. Since the vertices of $W$ are pairwise adjacent, the non-empty cliques of finite vertices are exactly the non-empty subsets $C\subseteq W$. If $|C|=j$, then
\[
    \lk(C,\flag{\Gamma_k})=\ast_{w_i\notin C}\{w_i,a_i\}\cong S^{k-j-1},
\]
and $\lk(C,\flag{\Gamma_k})_\infty$ is the simplex spanned by $\{a_i: w_i\notin C\}$. In particular, both complexes are empty when $C=W$. Moreover, $\St_{\Gamma_k}(C)$ is the join of the clique $C$ with a copy of $\Gamma_{k-j}$ (spanned by the pairs $\{w_i,a_i\}$ with $w_i\notin C$), where $\Gamma_0$ is the empty graph.

\begin{proposition}\label[proposition]{prop:products-virtually-free}
Let $R$ be a nonzero commutative ring and let $k\geq1$.
\begin{enumerate}
    \item $\BB_{\Gamma_k}$ is strictly of type $\fp_{k-1}(R)$, and it is of type $\textup{F}_{k-1}$.
    \item Let $1\neq Q\leqslant\BB_{\Gamma_k}$ be finite, and choose $g\in G(\Gamma_k)$ with $gQg^{-1}\leqslant[W]$, and put $\calC=\supp_{\Gamma_k}(gQg^{-1})\subseteq W$ and $j=|\calC|$. If $j<k$, then $N_{\BB_{\Gamma_k}}Q$ is strictly of type $\fp_{k-j-1}(R)$ and of type $\textup{F}_{k-j-1}$. If $j=k$, then $N_{\BB_{\Gamma_k}}Q\cong(\mathbb{Z}/2\mathbb{Z})^k$ is finite.
    \item The $G(\Gamma_k)$-conjugacy class of a finite subgroup $1\neq Q\leqslant\BB_{\Gamma_k}$ is a single $\BB_{\Gamma_k}$-conjugacy class if $j<k$, and splits into infinitely many $\BB_{\Gamma_k}$-conjugacy classes if $j=k$. In particular, $\BB_{\Gamma_k}$ has infinitely many conjugacy classes of finite subgroups and is not of type $\underline{\fp}_0(R)$.
\end{enumerate}
\end{proposition}

\begin{proof}
(1) In $\langle s_i\rangle*\langle t_i\rangle$, the kernel of the homomorphism onto $\mathbb{Z}/2\mathbb{Z}$ that sends $s_i$ to the generator and $t_i$ to $0$ has index $2$ and is free on $x_i=t_i$ and $y_i=s_it_is_i$. Hence
\[
    P=\prod_{i=1}^k\langle x_i,y_i\rangle\cong F_2^k
\]
is a torsion-free subgroup of index $2^k$ in $G(\Gamma_k)$. It is the right-angled Artin group $A_{O_k}$ on the flag complex $O_k=\ast_{i=1}^k\{x_i,y_i\}\cong S^{k-1}$, and since $f(x_i)=f(y_i)=1$, the restriction of $f$ to $P$ is the Bestvina--Brady homomorphism. Therefore $P\cap\BB_{\Gamma_k}=\BB_{O_k}$, which has finite index in $\BB_{\Gamma_k}$. The complex $O_k$ is $(k-2)$-connected and $\widetilde H_{k-1}(O_k;R)=R\neq0$, so by the classical Bestvina--Brady theorem the group $\BB_{O_k}$ is of type $\fp_{k-1}(R)$, but not of type $\fp_k(R)$; moreover, it is of type $\textup{F}_{k-1}$ by the homotopical form of that theorem \cite{BB}, which asserts that $\BB_L$ is of type $\textup{F}_{m+1}$ whenever $L$ is $m$-connected. Since these properties are invariant under passing to finite-index subgroups and overgroups, (1) follows.

(2) By \cref{thm:normalizer:BB:revised}, $N_{\BB_{\Gamma_k}}Q\cong\BB_{\St_{\Gamma_k}(\calC)}$. Since $\St_{\Gamma_k}(\calC)$ is the join of $\calC$ with a copy of $\Gamma_{k-j}$, we have $G(\St_{\Gamma_k}(\calC))=[\calC]\times G(\Gamma_{k-j})$, and $f$ vanishes on $[\calC]\cong(\mathbb{Z}/2\mathbb{Z})^j$. Hence
\[
    N_{\BB_{\Gamma_k}}Q\cong(\mathbb{Z}/2\mathbb{Z})^j\times\BB_{\Gamma_{k-j}}.
\]
If $j=k$, the second factor is trivial. If $j<k$, then $\BB_{\Gamma_{k-j}}$ has finite index in this group and the claim follows from (1) applied to $k-j$.

(3) By the computation above, $\lk(C,\flag{\Gamma_k})_\infty\neq\varnothing$ if and only if $C\neq W$. The first assertion then follows from the proof of \cref{thm:finite-conjugacy-classes-BB}; for instance, the subgroups $t_1^m[W]t_1^{-m}$, $m\in\mathbb{Z}$, are pairwise non-conjugate in $\BB_{\Gamma_k}$. The last assertion follows from \cref{thm:finite-conjugacy-classes-BB,thm:bredon-fp-criterion}.
\end{proof}

For instance, for $k=3$ the group $\BB_{\Gamma_3}$ is finitely presented and strictly of type $\fp_2(R)$, the normalizers of $\langle s_1\rangle$ and $\langle s_1,s_2\rangle$ are strictly of type $\fp_1(R)$ and $\fp_0(R)$, respectively, and that of $\langle s_1,s_2,s_3\rangle$ is finite.

\section*{AI declaration}
Most of the project was developed before the AI era. OpenAI Prism with GPT-5.2, along with Claude Opus 5.5, through a locally developed software tool, was used to assist in the development of the final example, to improve the editing and formatting of the manuscript, and to perform consistency checks throughout the text. The authors reviewed all AI-assisted suggestions incorporated into the manuscript and take full responsibility for its contents.

\bibliographystyle{alpha}
\bibliography{myblib}
\end{document}